\documentclass[a4paper,11pt, twoside]{article}

\usepackage[english]{babel}
\usepackage[latin1]{inputenc}
\usepackage[T1]{fontenc}
\usepackage{amsthm,amsmath,amssymb,amsfonts,graphics,graphicx}
\usepackage{color}
\usepackage{longtable}
\usepackage{natbib}
\def\var{\mathrm{Var}}

\def\argmin{\mathrm{argmin}}
\newcommand{\normp}[1]{\ensuremath{\vert\!\vert #1 \vert\!\vert_{2}}}
\newcommand{\indic}{1}
\newtheorem{Th}{Theorem}
\newtheorem{Remark}{Remark}
\newtheorem{Def}{Definition}
\newtheorem{Prop}{Proposition}
\newtheorem{Lemma}{Lemma}

\renewenvironment{proof}{\noindent{\bf Proof.}}{\hfill
  $\blacksquare$\par\noindent}
 \newcommand{\com}[1]{}
\newcommand{\norm}[1]{\ensuremath{\vert\!\vert #1 \vert\!\vert}}
\newcommand{\E}{\ensuremath{\mathbb{E}}}

\renewcommand{\P}{\ensuremath{\mathbb{P}}}
\newcommand{\X}{\ensuremath{\mathbb{R}}}
\newcommand{\R}{\ensuremath{\mathbb{R}}}

\newcommand{\Z}{\ensuremath{\mathbb{Z}}}

\newcommand{\Xx}{\ensuremath{\mathcal{X}}}

\renewcommand{\ln}{{\log\,}}
\renewcommand{\L}{\ensuremath{\mathbb{L}}}
\newcommand{\La}{\ensuremath{\Lambda}}
\newcommand{\Ga}{\ensuremath{\Gamma}}
\newcommand{\ga}{\ensuremath{\gamma}}
\newcommand{\al}{\ensuremath{\alpha}}
\newcommand{\la}{\ensuremath{\lambda}}

\newcommand{\e}{\ensuremath{\varepsilon}}
\newcommand{\p}{\ensuremath{\varphi}} 
\newcommand{\tp}{\ensuremath{{\tilde{\varphi}}}}
\newcommand{\be}{\ensuremath{\beta}}
\newcommand{\tb}{\ensuremath{\tilde{\beta}}}
\newcommand{\hb}{\ensuremath{\hat{\beta}}}
\newcommand{\jk}{\ensuremath{{j,k}}}
\newcommand{\N}[1]{\ensuremath{N_{#1}}}

\newcommand{\supp}{\ensuremath{\mbox{supp}}}

\renewcommand{\com}[1]{}
\numberwithin{equation}{section}

\newcommand{\pa}[1]{\textcolor{black}{#1}}

\title{Calibration of thresholding rules for Poisson  intensity estimation}
\date{}
\begin{document} 
%\maketitle

\noindent{\bf  \LARGE  Calibration of thresholding rules for Poisson}
\thispagestyle{empty}

\vspace{0.2cm}

\noindent{\bf \LARGE intensity estimation}

\vspace{0.6cm}

\noindent{\bf   \Large  Patricia  Reynaud-Bouret\footnotemark[1]  and  Vincent
Rivoirard\footnotemark[2]\footnotetext[1]{Laboratoire J. A. Dieudonn\'e, CNRS UMR 6621,
Universit\'e de Nice Sophia-Antipolis, Parc Valrose, 06108 Nice Cedex 2,
 France. Email:    reynaud@dma.ens.fr}\footnotetext[2]{Laboratoire de Math\'ematique,  CNRS UMR 8628,
Universit\'e Paris Sud, 91405 Orsay Cedex, France. D\'epartement
de Math\'ematiques et Applications, ENS-Paris, 45 Rue d'Ulm, 75230 Paris Cedex
05, France. Email: Vincent.Rivoirard@math.u-psud.fr}}

\vspace{0.8cm}

\noindent 

\begin{abstract}
In this paper, we deal with the problem of calibrating thresholding rules in the setting of Poisson intensity estimation. By using sharp concentration inequalities, oracle inequalities are derived and \pa{we }establish the optimality of our estimate up to a logarithmic term. This result is proved under mild assumptions and we do not impose any condition on the support of the signal to be estimated. Our procedure is based on data-driven thresholds. As usual, they depend on a threshold parameter $\gamma$ whose optimal value is hard to estimate from the data. Our main concern is to provide some theoretical and numerical results to handle this issue. In particular, we establish the existence of a minimal threshold parameter from the theoretical point of view: taking $\gamma<1$ deteriorates oracle performances of our procedure. In the same spirit, we establish the existence of a maximal threshold parameter and our theoretical results point out the optimal range $\gamma\in[1,12]$. Then, we 
 lead a numerical study that shows that choosing $\gamma$ larger than 1 but close to 1 is a fairly good choice. Finally, we compare our procedure with classical ones revealing the harmful role of the support of functions when estimated by classical procedures.
\end{abstract}

\noindent {\bf Keywords } Adaptive estimation, Calibration, Oracle inequalities, Poisson process, Wavelet thresholding \\

\noindent {\bf Mathematics Subject Classification (2000)} 62G05 62G20 
%%%%%%%%%%%%%%%%%%%%%%%%%%%%%%%%%%%%%%%%%%%
%%%%%%%%%%%%%%%%%%%%%%%%%%%%%%%%%%%%%%%%%%%
\section{Introduction}
In this paper, we consider the problem of estimating the intensity of a Poisson process. From a practical point of view, various methodologies have already been proposed. See for instance Rudemo \cite{rud} who proposed kernel and data-driven histogram rules calibrated by cross-validation. Thresholding algorithms have been performed by Donoho \cite{don} who modified the universal thresholding procedure proposed in \cite{dojo} by using the Anscombe transform or by Kolaczyk \cite{kol} whose procedure is based on the tails of the distribution of the noisy wavelet coefficients of the intensity. Finally, let us cite \pa{penalized} model selection type estimators built by Willett and Nowak \cite{wn} based on models spanned by piecewise polynomials. From the theoretical point of view, Cavalier and Koo \cite{ck} derived minimax rates on Besov balls by using wavelet thresholding. In the oracle approach, various optimal adaptive model selection rules have also been built by Baraud and Birg\'
 e \cite{bb}, Birg\'e \cite{bir} and Reynaud-Bouret \cite{ptrfpois}. Let us mention that these procedures are also minimax provided the intensity to be estimated is assumed to be supported by $[0,1]$. 

In a previous paper, \pa{we}  refined classical wavelet thresholding algorithms by proposing local data-driven thresholds (see \cite{Poisson_minimax}). Under very mild assumptions, \pa{the corresponding} procedure achieves optimal oracle inequalities and optimal minimax rates up to a logarithmic term. In particular, \pa{these} results are true even if the support of the intensity is unknown or infinite, which is rarely considered in the literature.  In \cite{Poisson_minimax}, we give many arguments to justify this unusual setting and we illustrate the influence of the support on minimax rates by showing how these rates deteriorate when the sparsity of the intensity decreases.  So, \pa{this}  algorithm, that is easily implementable, automatically adapts to the unknown regularity of the signal as usual, but also to the unknown support which is not classical. The main goal of this paper is to study the optimal calibration of the procedure studied in \cite{Poisson_minimax} from both theoretical and practical points of view. For this purpose, the next subsection briefly describes this procedure (Section \ref{data-driven} gives accurate definitions) and Section \ref{calissue} presents the calibration issue.
%%%%%%%%%%%%%%%%%%%%%%%%%%%%%%%%%%%%%%%%
\subsection{A brief description of our procedure}\label{brifdesc}
We observe a Poisson process $N$ whose mean measure $\mu$ is finite on the real line $\R$ and is  absolutely continuous with respect to the Lebesgue
measure (see  Section \ref{rappelPoisson} where  we recall classical  facts on
Poisson processes). Given $n$ a positive integer, we define  the intensity
of $N$ as the function $f$ that satisfies
$$f(x)=\frac{d\mu_x}{ndx}.$$
So, the total number of points of the process $N$, denoted $\mbox{card}(N)$, satisfies $$\E[\mbox{card}(N)]=n\norm{f}_1<\infty.$$ In particular, $\mbox{card}(N)$ is finite almost surely.  In the sequel, $f$ will be held
fixed and $n$  will go to $+\infty$. 
The
introduction of $n$  could seem artificial, but it allows to  present the following asymptotic theoretical  results in a  meaningful way since the mean of the number of points of $N$ goes to $\infty$ when $n\to\infty$.  In addition, our
framework is equivalent to the
observation of  a $n$-sample  of a Poisson  process with common  intensity $f$
with respect to  the Lebesgue measure.  The goal of this paper is to estimate $f$ by \pa{observing}  the points of $N$. 

First, we decompose the signal $f$ to be estimated as follows:
$$f=\sum_{\la\in \La}\be_\la \tp_\la \quad\mbox{  with }\quad \be_\la=\int \p_\la(x) f(x)
dx,$$
where  $((\p_\la)_{\la\in \La},(\tp_\la)_{\la\in \La})$ denotes a biorthogonal wavelet basis. In our paper, we mainly focus on the Haar basis (in this case, $\tp_\la=\p_\la$ for any $\la$) or on a special case of biorthogonal spline wavelet bases (in this case,  $\p_\la$ is piecewise constant and $\tp_\la$ is regular). See Section \ref{biorthogonal} where we recall well-known facts on biorthogonal wavelet bases or Cohen, Daubechies and Feauveau \cite{cdf} for a complete overview on such families. As usual in the wavelet setting, our goal is to estimate the wavelet coefficients $(\be_\la)_{\la}$ by thresholding empirical wavelet coefficients $(\hb_\la)_{\la}$ defined as
$$\hb_\la=\frac{1}{n}\sum_{T\in N}\p_\la(T).$$ 
Thresholding procedures have been introduced by Donoho and Johnstone \cite{dojo}.
Their main idea is that it is sufficient to keep a small amount
of the coefficients  to have a good estimation of the
function $f$. In our setting, the estimate of $f$ takes the form
$$\tilde{f}_{n,\gamma}=\sum_{\la \in \Gamma_n} \hb_\la\indic_{ \{|\hb_\la|\geq\eta_{\la,\gamma}\}}\tp_\la,$$
where $\Gamma_n$ is defined in (\ref{Gamman}). The thresholding procedure is detailed and discussed in Section \ref{data-driven}. We just mention here the form of the data-driven threshold $\eta_{\la,\ga}$:
$$\eta_{\la,\ga}=\sqrt{2\gamma \tilde{V}_{\la,n}\ln n }+\frac{\gamma\ln
    n}{3n}\norm{\p_\la}_\infty,$$
where  $\tilde{V}_{\la,n}$ is a sharp estimate of $\var(\hb_\la)$ defined in (\ref{Vtilde}) and \pa{where} $\ga$ \pa{is} a constant to be chosen. As explained in Section \ref{data-driven}, we have for most of the indices $\la$'s playing a key role for estimation:
$$\eta_{\la,\gamma}\approx\sqrt{2\gamma \tilde V_{\la,n}\ln n }.$$
In this case, $\eta_{\la,\ga}$ has a form close to the universal threshold $\eta^U$ proposed by Donoho and Johnstone \cite{dojo} in the Gaussian regression  framework:
$$\eta^U=\sqrt{2\sigma^2\log{n}},$$ where  $\sigma^2$  (assumed
to be known \pa{in the Gaussian framework}) is the variance of each noisy wavelet coefficient. Note, however, that our procedure depends on the so-called \emph{threshold parameter} $\ga$ that has to be properly chosen. The next section which describes calibration issues in a general way discusses this question.
%%%%%%%%%%%%%%%%%%%%%%%%%%%%%%%%%%%%
\subsection{The calibration issue}\label{calissue}
The major concern of this paper is the study of the calibration of the threshold parameter $\ga$: how should this parameter be chosen to obtain good results in both theory and practice? As usual, it can be proved that $\tilde{f}_{n,\gamma}$ achieves good theoretical performances in minimax or oracle points of view (see \cite{Poisson_minimax} or Theorem  \ref{inegoraclelavraie}) provided $\ga$ is large enough.  Such an assumption is very classical in the literature (see for
instance   \cite{aut}, \cite{ck},  \cite{djkp} or  \cite{jll}). Unfortunately, most of the time, the theoretical choice of the threshold parameter is not suitable for practical issues. More precisely, this choice is often too conservative.  See for instance Juditsky and
Lambert-Lacroix \cite{jll} who illustrate this statement in Remark 5 of their paper: their threshold  parameter, denoted $\la$, has to be larger than 14 to obtain theoretical results, but they suggest to use $\la\in [\sqrt{2},2]$ for practical issues. So, one of the main goals of this paper is to fill the gap between the optimal parameter choice provided by theoretical results on the one hand and by a simulation study on the other hand.

Only a few papers have been devoted to theoretical calibration of statistical procedures. In the model selection setting, the issue of calibration has been addressed by Birg\'e and Massart \cite{bm}. They considered  penalized estimators in a Gaussian homoscedastic regression framework with known variance and calibration of penalty constants is based on the following methodology.  They showed that there exists a minimal penalty in the sense that taking smaller penalties leads to inconsistent estimation procedures. Under some conditions, they further prove that the optimal penalty is twice the minimal penalty. This relationship characterizes the ``slope heuristic'' of Birg\'e and Massart \cite{bm}. Such a method has been successfully applied for practical purposes in \cite{leb}. Baraud, Giraud and Huet \cite{bgh} (respectively Arlot and Massart \cite{am}) generalized these results when the variance is unknown (respectively for non-Gaussian or heteroscedastic data).
These approaches constitute alternatives to popular cross-validation methods (see \cite{all} or \cite{sto}). For instance, $V$-fold cross-validation (see \cite{gei}) is widely used to calibrate procedure parameters but its computational cost can be high.

%%%%%%%%%%%%%%%%%%%%%%%%%%%%%%%%%%%%
\subsection{Our results}
The starting point of our results is the oracle inequality stated in Section \ref{data-driven}:  Theorem \ref{inegoraclelavraie} shows that the estimate  $\tilde{f}_{n,\gamma}$ achieves the oracle risk up to a logarithmic term. This result is true as soon as $\gamma>1$ and $f\in\L_2\cap\L_1$. In particular, nothing is assumed with respect to the support of $f$ or $\norm{f}_\infty$: 
our result remains true if $\norm{f}_\infty=\infty$ and if the support of $f$ is unknown or infinite. The oracle inequality of Theorem \ref{inegoraclelavraie} is refined in Section \ref{study} where $f$ is assumed to belong to a special class denoted ${\cal F}_n(R)$ whose signals have only a finite number of non-zero wavelet coefficients (see Theorem \ref{classFn}).

Then, in the perspective of calibrating thresholding rules, we consider theoretical performances of $\tilde{f}_{n,\gamma}$ with $\ga<1$ by using the Haar basis. For the signal $f=\indic_{[0,1]}$, Theorem \ref{inegoraclelavraie} shows that $\tilde{f}_{n,\gamma}$ with $\ga>1$ achieves the rate $\frac{\log n}{n}$. But the lower bound of Theorem \ref{lower} shows that the rate of $\tilde{f}_{n,\gamma}$ with $\ga<1$ is larger than $n^{-\delta}$ for $\delta<1$. So, as in \cite{bm} for instance, we prove  the existence of a minimal threshold parameter: $\gamma=1$. Of course, the next step concerns the existence of a maximal threshold parameter. This issue is answered by Theorem \ref{uppth} which studies the maximal ratio between the risk of $\tilde{f}_{n,\gamma}$ and the oracle risk on ${\cal F}_n(R)$. We derive a lower bound that shows that taking $\ga>12$ leads to worse rates constants: \pa{this is consequently} a bad choice.

The optimal choice for $\gamma$ is derived from a numerical study, keeping in mind that the theory points out the range $\gamma\in[1,12]$. Some simulations are provided for estimating various signals by considering either the Haar basis or a particular biorthogonal spline wavelet basis (see Section \ref{simulations}). Our numerical results show that choosing $\gamma$ larger than 1 but close to 1 is a fairly good choice, which corroborates theoretical results. Actually, our simulation study suggests that Theorem \ref{lower} remains true for all signals of ${\cal F}_n(R)$ whatever the basis for decomposing signals is used.% (the proof of this conjecture constitutes futher works).

Finally, we lead a comparative study with other competitive procedures. We show that the thresholding rule proposed in this paper outperforms universal thresholding (when combined with the Anscombe transform) or Kolaczyk's procedure. Finally, the robustness of our procedure with respect to the support issue is emphasized and we show the harmful role played by large supports of signals when estimation is performed by other classical procedures.

%%%%%%%%%%%%%%%%%%%%%%%%%%%%%%%%%%%%
\subsection{Overview of the paper}
Section \ref{data-driven} defines the thresholding estimate $\tilde{f}_{n,\gamma}$ and studies its properties under the oracle approach. In Section \ref{study}, we refine this study on the set of positive functions that can be decomposed on a finite combination of the basis. Calibration of thresholds is discussed in Section \ref{penaltyterm} and Section \ref{simulations} illustrates
our theoretical results by some simulations. Section \ref{proofs} is devoted to the proofs of the results. Finally, Section \ref{maintools} recalls well-known facts on Poisson processes and biorthogonal wavelet bases.

\section{Data-driven thresholding rules and oracle inequalities}\label{data-driven}
The goal of this section is to specify our thresholding rule. For this  purpose, we  assume that $f$  belongs to  $\L_2(\R)$ and we  use the
decomposition of $f$ on one of the biorthogonal wavelet bases described in Section \ref{biorthogonal}.  We recall that, as classical
orthonormal  wavelet  bases,  biorthogonal  wavelet  bases  are  generated  by
dilatations  and translations  of father  and mother  wavelets.  But considering
biorthogonal  wavelets  allows  to  distinguish, if  necessary,  wavelets  for
analysis (that are piecewise constant functions in this paper) and wavelets for
reconstruction     with     a      prescribed     number     of     continuous
derivatives. Then, the decomposition of $f$ on a biorthogonal wavelet basis takes
the following form:
\begin{equation}\label{decom2}
f=\sum_{k\in\Z}\alpha_k\tilde\phi_k+\sum_{j\geq
0}\sum_{k\in\Z}\beta_{j,k}\tilde\psi_{j,k},
\end{equation}
where for any  $j\geq 0$ and any $k\in\Z$,
$$\alpha_k=\int_\R          f(x)\phi_k(x)dx,\quad          \beta_{j,k}=\int_\R
f(x)\psi_{j,k}(x)dx.$$
 See Section \ref{biorthogonal} for further details. To shorten mathematical
expressions, we set
\begin{equation}\label{Lajk}
\Lambda=\{\la=(j,k):\quad j\geq -1,k\in\Z\}\nonumber
\end{equation}
and  for   any  $\la\in\Lambda$,   $\p_\la=\phi_k$  (respectively
$\tilde\p_\la=\tilde\phi_k$) if $\la=(-1,k)$ and
$\p_\la=\psi_{j,k}$ (respectively $\tilde\p_\la=\tilde\psi_{j,k}$) if $\la=(j,k)$ with $j\geq 0$. Similarly, $\be_\la=\al_k$
if   $\la=(-1,k)$   and  $\be_\la=\be_{j,k}$   if   $\la=(j,k)$  with   $j\geq
0$. Now, (\ref{decom2}) can be rewritten as
\begin{equation}
\label{decom}
f=\sum_{\la\in \La}\be_\la \tp_\la \quad\mbox{  with }\quad \be_\la=\int \p_\la(x) f(x)
dx.
\end{equation}
In particular,  (\ref{decom}) holds for  the Haar basis that will play a special role in this paper,  where in this
case $\tilde\p_\la=\p_\la$.
Now, let us define the thresholding  estimate of $f$ by using the properties of
Poisson processes.  First, we introduce for any $\la\in\Lambda$, the natural estimator of $\be_\la$ defined by
\begin{equation}
\label{defest1}
\hb_\la=\frac{1}{n}\int \p_\la(x) dN_x,
\end{equation}
where   we  denote  by   $dN$  the   discrete  random   measure  $\sum_{T\in N}\delta_T$ and  for any compactly supported function $g$,
$$\int  g(x)  dN_x=\sum_{T\in  N}g(T).$$
So, the estimator $\hb_\la$ is unbiased: $\E(\hb_\la)=\be_\la$. Then, given some parameter $\ga>0$, we define the threshold
$\eta_{\la,\ga}$ mentioned in Introduction as
\begin{equation}\label{defthresh}
\eta_{\la,\ga}=\sqrt{2\gamma \tilde{V}_{\la,n}\ln n }+\frac{\gamma\ln
    n}{3n}\norm{\p_\la}_\infty,
\end{equation}
with 
\begin{equation}\label{Vtilde}
\tilde{V}_{\la,n}=\hat{V}_{\la,n}+\sqrt{2\gamma     \ln    n    \hat{V}_{\la,n}
\frac{\norm{\p_\la}_\infty^2}{n^2}}+3                \gamma                \ln
n\frac{\norm{\p_\la}_\infty^2}{n^2}
\end{equation}
where $$\hat{V}_{\la,n}=\frac{1}{n^2}\int \p_\la^2(x) dN_x.$$
Note that $\hat{V}_{\la,n}$ satisfies $\E(\hat{V}_{\la,n})=V_{\la,n}$,
where $$V_{\la,n}=\var(\hb_\la)=\frac{1}{n} \int \p_\la^2(x) f(x) dx.$$
Finally, with
\begin{equation}\label{Gamman}
\Gamma_n=\left\{\la=(j,k)\in\Lambda:\quad j\leq j_0\right\}
,
\end{equation}
 where $j_0=j_0(n)$ is the   integer   such   that  $2^{j_0}\leq n<2^{j_0+1}$,
we set for any $\la\in\Lambda$,
$$\tb_\la=\hb_\la\indic_{\{|\hb_\la|\geq\eta_{\la,\gamma}\}}\indic_{\{\la\in\Gamma_n\}}$$
and $\tb=(\tb_\la)_{\la\in\Lambda}$.  Finally, the estimator of $f$ is
\begin{equation}\label{defest}
\tilde{f}_{n,\gamma}=\sum_{\la \in \La} \tb_\la \tp_\la
\end{equation}
and only depends on the choice of $\gamma$. 
  When the Haar basis is used, the estimate is
denoted $\tilde{f}_{n,\gamma}^H$ and its wavelet coefficients are denoted $\tb^H=(\tb_\la^H)_{\la\in\Lambda}$.  The threshold $\eta_{\la,\gamma}$ seems to be defined in a rather
complicated manner but we can notice the following fact. Given $\lambda\in\Gamma_n$, when  there exists  a constant
$c_0>0$ such that $f(x)\geq
c_0$ for $x$ in the support of $\p_\la$ \pa{satisfying} %and %if $\lambda\in\Gamma_n$ satisfies
 $\|\p_\la\|_{\infty}^2=o_n(n(\log   n)^{-1})$, \pa{then},  with large probability, the   deterministic  term   of
(\ref{defthresh}) is negligible with respect to the random one. \pa{In this case} we asymptotically \pa{derive}
\begin{equation}\label{approxseuil}
\eta_{\la,\gamma}\approx\sqrt{2\gamma \tilde V_{\la,n}\ln n },
\end{equation}
as stated in Introduction. 
Actually, the deterministic term of
(\ref{defthresh}) allows to consider $\gamma$ close
to 1 and to control large deviations terms for high resolution levels. In the same spirit, $V_{\la,n}$ is slightly
overestimated and  we consider $\tilde V_{\la,n}$ instead  of $\hat V_{\la,n}$
to   define  the  threshold. 

The performance of this procedure has been investigated in the oracle point of view in \cite{Poisson_minimax}. We recall that in  the context of  wavelet function
estimation by thresholding, the oracle does not tell us the true function, but 
tells us  the coefficients that have to be kept. This ``estimator'' obtained
with the aid of an oracle is not a true estimator, of course, since it depends
on $f$.  But it represents  an ideal for the particular estimation  method. The goal of the oracle approach  is   to  derive  true  estimators  which  can essentially ``mimic'' the performance of the ``oracle estimator''. 
 In our framework, it is easy to see  that the oracle estimate is $\bar{f}=\sum_{\la\in \Ga_n}
\bar{\be}_\la\tp_\la$, where %for any $\la\in \Ga_n$,
$\bar{\be}_\la=\hb_\la \indic_{\{\be_\la^2>V_{\la,n}\}}$
%and we have
\pa{satisfies}
$$\E((\bar{\be}_\la-{\be}_\la)^2)=\min(\be_\la^2,V_{\la,n}).$$
By keeping the coefficients $\hb_\la$ larger than \pa{the} thresholds defined in (\ref{defthresh}),  our estimator has a risk that is not larger than the oracle risk,  up to a logarithmic term, as stated by the following key result.
\begin{Th}
\label{inegoraclelavraie} 
Let us consider a biorthogonal wavelet basis satisfying the properties described in
Section \ref{biorthogonal}.   If $\ga>1$,  then $\tilde  f_{n,\gamma}$ satisfies  the
following oracle inequality: for $n$ large enough
\begin{equation}\label{inegoraclelavraie1}
\E(\norm{\tilde{f}_{n,\gamma}-f}_2^2)\leq
C_1\ln n\sum_{\la\in\Gamma_n}\min(\be_\la^2,V_{\la,n})+C_1\sum_{\la\notin\Gamma_n}\be_\la^2+\frac{C_2}{n}
\end{equation}
where $C_1$  is a positive  constant depending  only on
$\ga$ and on the functions that generate the biorthogonal wavelet basis. $C_2$ is also a
positive constant depending on
$\ga$, $\|f\|_1 $ and on the functions that generate the basis. 
\end{Th}
 Following the oracle point of view of Donoho and Johnstone, Theorem
\ref{inegoraclelavraie} shows that our procedure is optimal up to the logarithmic factor. This logarithmic term is in some sense unavoidable. It is the price we pay for adaptivity (i.e. for not  knowing the coefficients that we must keep). Our result is true provided $f\in\L_1(\R)\cap\L_2(\R)$. So, assumptions on $f$ are very mild \pa{here. This is not the case for} %Whereas 
most of the results for non-parametric estimation procedures \pa{where} one assumes that $\norm{f}_\infty<\infty$ and \pa{that} $f$ has a compact support. \pa{Note in addition that} this support and $\norm{f}_\infty$ are often known in the literature. %, we do not need these conditions to state the oracle result. In particular, 
 \pa{On the contrary, in Theorem \ref{inegoraclelavraie}} $f$ and its support can be unbounded. So, we make as few assumptions as possible. This is allowed by considering random thresholding with the data-driven thresholds defined in (\ref{defthresh}). This result is proved in \cite{Poisson_minimax} where in addition optimality properties of the estimate (\ref{defest}) under the minimax approach are established.

%Now  that the  constants are  precisely written,  we may  ask ourself  how to
%choose $\ga$ such that the upper bound given by Theorem \ref{inegmodelsel} is
%the  sharpest  for  $\tilde{f}_\ga$.  The  problem  is  that  there  are  two
%quantities in  the threshold $\eta_{\la,\ga}$ (\ref{defthresh})  and that the
%last term, which can be viewed in a first approach as negligible with respect
%to the  main square root term, is  precisely here to deal  in particular with
%the cases where the mass $F_\la$ is too small.  In order to remove the cases where the mass becomes too small, we restrict ourselves to
A glance at the proof of Theorem
\ref{inegoraclelavraie} shows that the constants $C_1$ and $C_2$ strongly depends on $\gamma$. Actually, without further assumptions on $f$, the constants $C_1$ and $C_2$ blow up when $\gamma$ tends to $1$. In particular, such an oracle inequality is not sharp enough for some calibration issues. In the next section, we investigate this problem and we derive sharp oracle inequalities for a large class of functions. Furthermore, the upper bound in (\ref{majoFn})   depends on absolute constants whose size is acceptable.
%%%%%%%%%%%%%%%%%%%%%%%%%%%%%%%%%%%%%%%%%%%%%%%%%%%%%%
%%%%%%%%%%%%%%%%%%%%%%%%%%%%%%%%%%%%%%%%%%%%%%%%%%%%%%
\section{Study on a special class of functions}\label{study}
In the sequel, we consider the Haar basis and the estimator $\tilde{f}^H_{n,\gamma}$. We restrict our study on estimation of the functions of ${\cal F}$ defined as the set of positive functions that can be decomposed on a finite combination of $(\tilde\p_\la)_{\la\in\Lambda}$:
$${\cal F}=\left\{f=\sum_{\la\in \La}\be_\la \tp_\la\geq 0:\quad \mbox{card}\{\la\in \La: \ \be_\la\not=0\}<\infty\right\}.$$
To study sharp performances of our procedure, we introduce a subclass of the class ${\cal F}$: for any $n$ and any radius $R$, we define:
$${\cal F}_{n}(R)=\left\{f\geq 0:\quad f\in \L_1(R)\cap\L_2(R)\cap\L_{\infty}(R),\ F_\la\geq
\frac{(\ln            n)(\ln\ln           n)}{n}1_{\be_\la\not=0},           \
\forall\;\la\in\Lambda\right\},$$
where for any $\la$, we set $$F_\la=\int_{\supp(\p_\la)} f(x) dx\quad\mbox{and}\quad\supp(\p_\la)=\left\{x\in\R:\ \p_\la(x)\not=0\right\}, $$ which allows to establish a decomposition of ${\cal F}$. Indeed, we have the following result proved in Section \ref{preuvepropsharp}:
\begin{Prop}\label{Propsharp}
 When $n$ (or $R$) increases, $\left({\cal F}_{n}(R)\right)_{n,R}$ is a non-decreasing sequence of sets. In addition, we have:
$$\bigcup_{n}\bigcup_{R} {\mathcal F}_n(R)={\cal F}.$$
\end{Prop}
 The definition of ${\mathcal F}_n(R)$ especially relies on the technical condition
\begin{equation}\label{condFn}
F_\la\geq
\frac{(\ln            n)(\ln\ln           n)}{n}1_{\be_\la\not=0}.
\end{equation}
Remember that the distribution of the number of points of $N$ that lies in $\supp(\p_\la)$ is the Poisson distribution with mean $nF_\la$. So, the previous condition ensures that we have a significant number of %realizations 
\pa{points} of $N$ to estimate non-zero wavelet coefficients.
Another main point is that under (\ref{condFn}),
$$\sqrt{V_{\la,n}  \ln  n}  \geq  \frac{\ln  n
 \norm{\p_\la}_\infty}{n} \times\sqrt{\ln\ln n}$$
(see Section \ref{sec:proofclassFn}), so (\ref{approxseuil}) is true with large probability.
%which means that the  indices $\la$  for which  the
%wavelet coefficients $\be_\la$ of $f\in{\mathcal F}_n(R)$ are not equal to 0 correspond to the indices for which $F_\la$ is ''quite large'', so  $f$ is ''far''  from 0 on  $\supp(\p_\la)$.
%This allows to control the number of points of $N$ that lies in $\supp(\p_\la)$ and whose distribution is a Poisson distribution with mean $F_\la$. 
The term $\frac{(\ln            n)(\ln\ln           n)}{n}$ appears for technical reasons but could be replaced by any term $u_n$ such that
$$\lim_{n\to\infty}u_n=0\quad\mbox{and}\quad\lim_{n\to\infty}u_n^{-1}\left(\frac{\ln n}{n}\right)=0.$$
 %and allows  to control the number of points of $N$ that lies in $\supp(\p_\la)$ and whose distribution is a Poisson distribution with mean equal to $F_\la$. More precisely, $F_\la$ has to be asymptotically larger than the term $\frac{\ln n}{n}$ that  is related to the rates established in Theorem \ref{classFn}. This is the reason  why we add the term $\ln\ln n$ that can be replaced by any term that goes to $\infty$ but slower than any power of $\log n$. Actually, the  indices $\la$  for which  the wavelet coefficients $\be_\la$ of $f\in{\mathcal F}_n(R)$ are not equal to 0 correspond to the indices for which $F_\la$ is ''quite large'', which means that  $f$ is ''far''  from 0 on  $\supp(\p_\la)$. 
In practice, many interesting signals are well approximated by a function of ${\mathcal F}$. 
So, using Proposition \ref{Propsharp},  a convenient estimate is an estimate
with a good behavior on ${\mathcal F}_n(R)$, at least for large values of $n$ and $R$. 
Furthermore, note that we do not have any restriction on the \pa{precise} location of the support of functions of ${\mathcal F}_n(R)$ (\pa{even if} %remember that 
these functions have only a finite set of non-zero wavelet coefficients). This provides a second reason for considering ${\mathcal F}_n(R)$ if we are interested in estimated signals with unknown or infinite supports. 
We  now focus on $\tilde{f}^H_{n,\gamma}$ with the special value $\ga =1+\sqrt{2}$ and we study its properties on ${\cal F}_n(R)$.
\begin{Th}\label{classFn}
Let $R>0$ be fixed. Let $\ga =1+\sqrt{2}$ and let $\eta_{\la,\gamma}$ be as in
(\ref{defthresh}). Then  $\tilde{f}^H_{n,\gamma}$ achieves the following oracle inequality: for  $n$ large enough, for any $f\in{\cal F}_{n}(R)$,
\begin{equation}\label{majoFn}
\E(\norm{\tilde{f}_{n,\gamma}^H-f}_2^2)\leq
12\ln n\left[\sum_{\la\in\Gamma_n}\min(\be_\la^2,V_{\la,n})+\frac{1}{n}\right].
%\sup_{f\in {\cal F}_n(R)}\frac{\E(\normp{\tilde{f}^H_{n,\gamma}-f}^2)}{\sum_{\la \in \Ga_n}\min(\be_\la^2, V_{\la,n})+\frac{1}{n}}\leq   12   \ln    n .
\end{equation}
\end{Th}
Inequality (\ref{majoFn}) shows that on ${\cal F}_n(R)$, our estimate achieves the oracle risk up to the term $12\log n$ and the \pa{negligible} term $\frac{1}{n}$. Finally, \pa{let us} mention that when $f\in{\cal F}_{n}(R)$,
$$\sum_{\la\notin\Gamma_n}\be_\la^2=0.$$
Our result is stated with $\ga=1+\sqrt{2}$. This value comes from optimizations of upper bounds given by Lemma \ref{lemmesharp} stated in Section \ref{sec:proofclassFn}. \pa{This} constitutes a first theoretical calibration result and \pa{this} is the first step for choosing the parameter $\ga$ in an optimal way. The next section further investigates this problem.
%%%%%%%%%%%%%%%%%%%%%%%%%%%%%%%%%%%%%%%%%%%%%%%%%%%%%%
%%%%%%%%%%%%%%%%%%%%%%%%%%%%%%%%%%%%%%%%%%%%%%%%%%%%%%
\section{How to choose the parameter $\gamma$}\label{penaltyterm}
%Now, we address the problem of choosing conveniently the threshold parameter $\gamma$ from the theoretical point of view.  
In this Section, our goal is to find lower and upper bounds for the parameter $\gamma$. 
%The aim and the proofs are inspired by \cite{bm} who considered penalized estimators and calibrated constants for penalties in a Gaussian framework. In particular, they showed that if the penalty constant is smaller than 1, then the penalized estimator behaves in a quite unsatisfactory way. This study was used in practice to derive adequate data-driven penalties by Lebarbier \cite{leb}. 
Theorem \ref{inegoraclelavraie} established that for any signal, we achieve the oracle estimator up to a logarithmic term provided $\gamma>1$. So, our primary interest is to wonder what happens, from the theoretical point of view, when $\ga\leq 1$? To handle this problem, we consider the simplest signal in our setting, namely
$$f=\indic_{[0,1]}.$$ Applying Theorem \ref{inegoraclelavraie} with the Haar basis and $\ga>1$ gives
$$\E(\normp{\tilde{f}_{n,\ga}^H-f}^2)\leq C\frac{\ln n}{n},$$
where $C$ is a constant. The following result shows that this rate cannot be achieved for this particular signal when $\ga<1$.
\begin{Th}\label{lower}
Let $f=\indic_{[0,1]}$.
If $\ga<1$ then there exists $\delta<1$ not dependent of $n$ such that 
$$\E(\norm{\tilde{f}_{n,\gamma}^H-f}^2_2)\geq \frac{c}{n^{\delta}},$$
where $c$ is a constant.
\end{Th}
Theorem \ref{lower}   establishes     that,    asymptotically,
$\tilde{f}^H_{n,\gamma}$       with      $\gamma<1$       cannot      estimate a very simple signal ($f=\indic_{[0,1]}$)  at a convenient rate  of convergence.  This  provides a
lower bound for the threshold parameter $\gamma$: we have to take  $\ga\geq 1$.

Now, let us study the upper bound for the parameter $\ga$. For this purpose, we do not consider a particular signal, but we use the worst oracle ratio on the whole class ${\cal F}_n(R)$. Remember that when $\ga=1+\sqrt{2}$, Theorem \ref{classFn} gives that this ratio cannot grow faster than $12 \ln    n$, when $n$ goes to $\infty$: for $n$ large enough,
$$\sup_{f\in {\cal F}_n(R)}\frac{\E(\normp{\tilde{f}^H_{n,\gamma}-f}^2)}{\sum_{\la \in \Ga_n}
\min(\be_\la^2, V_{\la,n})+\frac{1}{n}}\leq   12   \ln    n .$$  Our aim is to
establish that the oracle ratio on ${\cal F}_n(R)$ for the estimator $\tilde{f}^H_{n,\ga}$ where $\ga$ is large,  is larger than the previous upper bound. This
goal is reached in the following theorem.
 \begin{Th}\label{uppth}
Let $\ga_{\min}>1$ be fixed and let $\ga >\ga_{\min}$. 
Then, for any $R\pa{\geq 2}$,
$$\sup_{f\in {\cal F}_n(R)}\frac{\E(\normp{\tilde{f}^H_{n,\gamma}-f}^2)}{\sum_{\la \in \Ga_n}
\min( \beta_\la^2, V_{\la,n})+\frac{1}{n}}\geq 2(\sqrt{\ga}-\sqrt{\ga_{\min}})^2 \ln n \times(1+o_n(1)).$$
\end{Th}
Now, if we choose $\ga>(1+\sqrt{6})^2\approx 11.9$, we can take $\ga_{\min}>1$ such that the  resulting
maximal oracle ratio of $\tilde{f}^H_{n,\gamma}$
is larger than $12 \ln n$ for $n$ large enough. So, taking $\ga>12$ is a bad choice for estimation on the whole class ${\cal F}_n(R)$.

Note  that  the  function  $\indic_{[0,1]}$ belongs  to  $\mathcal{F}_n(\pa{2})$, \pa{for all  $n \geq 2$}.  So,
combining  Theorems  \ref{classFn}, \ref{lower} and \ref{uppth}
proves  that  the convenient  choice  for  $\gamma$  belongs to  the  interval
$[1,12]$. Finally, observe that  the rate exponent deteriorates for $\gamma<1$
whereas we only prove that the  choice $\gamma>12$ leads to worse
rates constants. 
%Let  us explain  more  precisely this  class.  Remember that  in our  problem
%$\la=(j,k)$. If $j$ is fixed and $k$  varies, one sees that if $f$ belongs to
%$\mathcal{F}_n$, then  $f$ has  a bounded support  depending on $n$  and that
%nothing  is said  on the  precise location  of the  support, it  may  tend to
%$\infty$ with $n$. Moreover $f$ is lower bounded on its support by a quantity
%becoming smaller with $n$. Finally as  $f$ is upper bounded it means that the
%level of resolution $2^j$ is bounded  by $n$. So from a model selection point
%of view it  means that $f$ is characterized by a  family of coefficients that
%are null outside  a particular $m$, which may change with  $f$, but still $m$
%remains a  finite subset  of $\{\la=(j,k), 2^j\leq  n\}$. So looking  at this
%problem from the other way around, if $f$ is characterized by a finite number
%of   coefficients,  then   for  $n$   large  enough,   $f$  will   belong  to
%$\mathcal{F}_n$ and the least that we  could ask from a model selection point
%of view is to  be able to estimate $f$ with the  best accuracy. Indeed if $f$
%belongs to one model, since the  models are thought to be nicer than anything
%else, a procedure that will not be as good as possible for this $f$ is in our
%opinion not a good one. 
%%%%%%%%%%%%%%%%%%%%%%%%%%%%%%%%%%%%%%%%%%%%%%%%%%%%%%
%%%%%%%%%%%%%%%%%%%%%%%%%%%%%%%%%%%%%%%%%%%%%%%%%%%%%%
\section{Numerical study}\label{simulations}
\begin{figure}[t]
\begin{center}
\includegraphics[width=0.7\linewidth,angle=0]{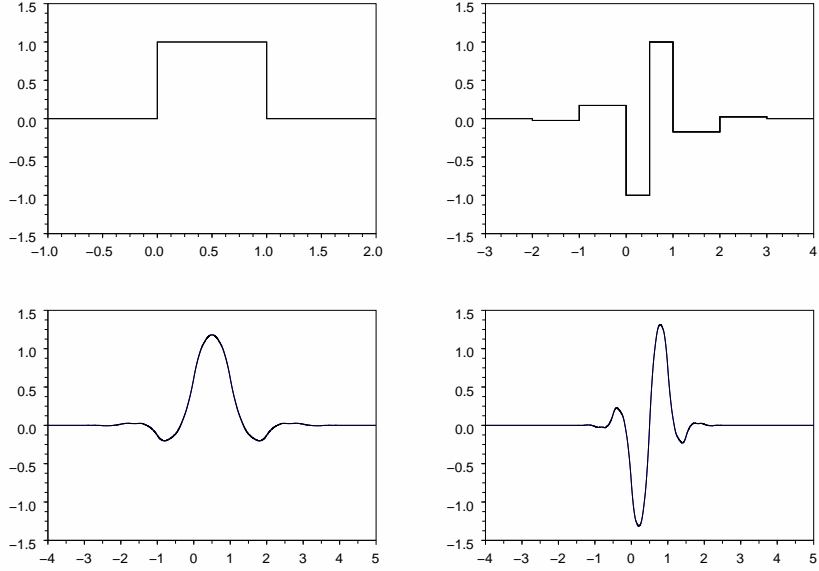}
\end{center}
\caption{The spline basis.  Top: $\phi$ and $\psi$, Bottom: $\tilde\phi$ and $\tilde\psi$}\label{fig-basespline}
\end{figure}

In this  section, some  simulations are provided  and the performances  of the
thresholding  rule are  measured from  the numerical  point of  view by comparing our estimator with other well-known procedures.  We also
discuss the ideal  choice for the parameter $\gamma$ keeping in mind that the
value $\gamma=1$ constitutes a border for the theoretical results (see Theorems
\ref{inegoraclelavraie} and \ref{lower}).  For these purposes, our procedure is performed for estimating various intensity
signals and the wavelet set-up associated with biorthogonal
wavelet bases is considered. More precisely, we focus either on the Haar basis where
$$\phi=\tilde\phi=\indic_{[0,1]},\quad
\psi=\tilde\psi=\indic_{[0,1/2]}-\indic_{]1/2,1]}$$or  on a  special  case of
spline systems given in Figure \ref{fig-basespline}. 
The latter, called hereafter the spline
basis,  has  the following  properties.  First,  the  support  of  $\phi$,  $\psi$,
$\tilde\phi$ and $\tilde\psi$ is included in $[-4,5]$. The reconstruction wavelets $\tilde\phi$ and $\tilde\psi$
belong to  $C^{1.272}$. Finally,  the wavelet  $\psi$ is  a piecewise
constant function orthogonal to polynomials of degree 4 (see \cite{don}). So, such
a basis has properties 1--5 required in Section 
\ref{biorthogonal} with $r=0.272$.  Then, the signal $f$ to be estimated is decomposed
as follows:
$$f=\sum_{\la\in\Lambda}\be_\la\tilde\p_\la=\sum_{k\in\Z}\be_{-1,k}\tilde\phi_{k}+\sum_{j\geq
0}\sum_{k\in\Z}\be_{j,k}\tilde\psi_{j,k}.$$
For  estimating $f$, we  use the empirical coefficients $(\hat\be_\la)_{\la\in\Lambda}$
associated with a Poisson process $N$ whose
intensity with respect to the Lebesgue measure is $n\times f$.
Since $\phi$ and $\psi$ are piecewise constant functions, accurate values of the \pa{empirical coefficients} are available, which allows to avoid many computational and approximation issues that often arise in the wavelet setting.
%To shed light on typical aspects of Poisson intensity estimation, Figure \ref{fig-observations} displays the reconstruction  obtained by using all the  noisy wavelet coefficients of  a particular signal  (the density of a Gaussian variable with mean 0.5 and standard deviation 0.25) with $n=4096$. We mean that $(\be_{j,k})_{j\geq -1,k\in\Z}$ is estimated by $(\hat\be_{j,k})_{-1\leq j\leq 10,k\in\Z}$ without using thresholding. 
%\begin{figure}[t]\begin{center}\includegraphics[width=0.7\linewidth,angle=0]{bruit.eps}\end{center}\caption{Plots of the signal $f(x)=\frac{1}{0.25\sqrt{2\pi}}\exp\left(\frac{(x-0.5)^2}{2\times 0.25^2}\right)$ and purely noisy reconstruction with $n=4096$ based on the wavelet coefficients until the level 10 and by using the Haar basis.}\label{fig-observations}\end{figure}
%As  expected, variability  highly depends  on the  local values of the signal.  So, our framework is very different from classical regression where we observe random variables with common variance.  
We consider the thresholding rule $\tilde      f_\gamma=(\tilde{f}_{n,\gamma})_n$      with
$\tilde{f}_{n,\gamma}$ defined in (\ref{defest}) with 
$$\Gamma_n=\left\{\la=(j,k):\quad -1\leq j\leq j_0,\ k\in\Z\right\}$$
and 
$$\eta_{\la,\gamma}=\sqrt{2\gamma\ln (n) \hat{V}_{\la,n} }+\frac{\gamma\ln
    n}{3n}\norm{\p_\la}_\infty.$$
Observe that  $\eta_{\la,\gamma}$ slightly differs from  the threshold defined
in    (\ref{defthresh})   since    $\tilde{V}_{\la,n}$    is now  replaced    with
$\hat{V}_{\la,n}$.  It allows to derive the parameter $\gamma$ as an explicit function of the
threshold which is necessary to draw figures without using a discretization of
$\gamma$, which is crucial in Section \ref{calibrationnumerical}. The performances of our
thresholding rule associated with the threshold $\eta_{\la,\gamma}$ defined in
(\ref{defthresh}) are probably equivalent (see (\ref{inegseuil})). 

\begin{figure}[tpb]
\begin{center}
\includegraphics[width=1.2\linewidth,angle=-90]{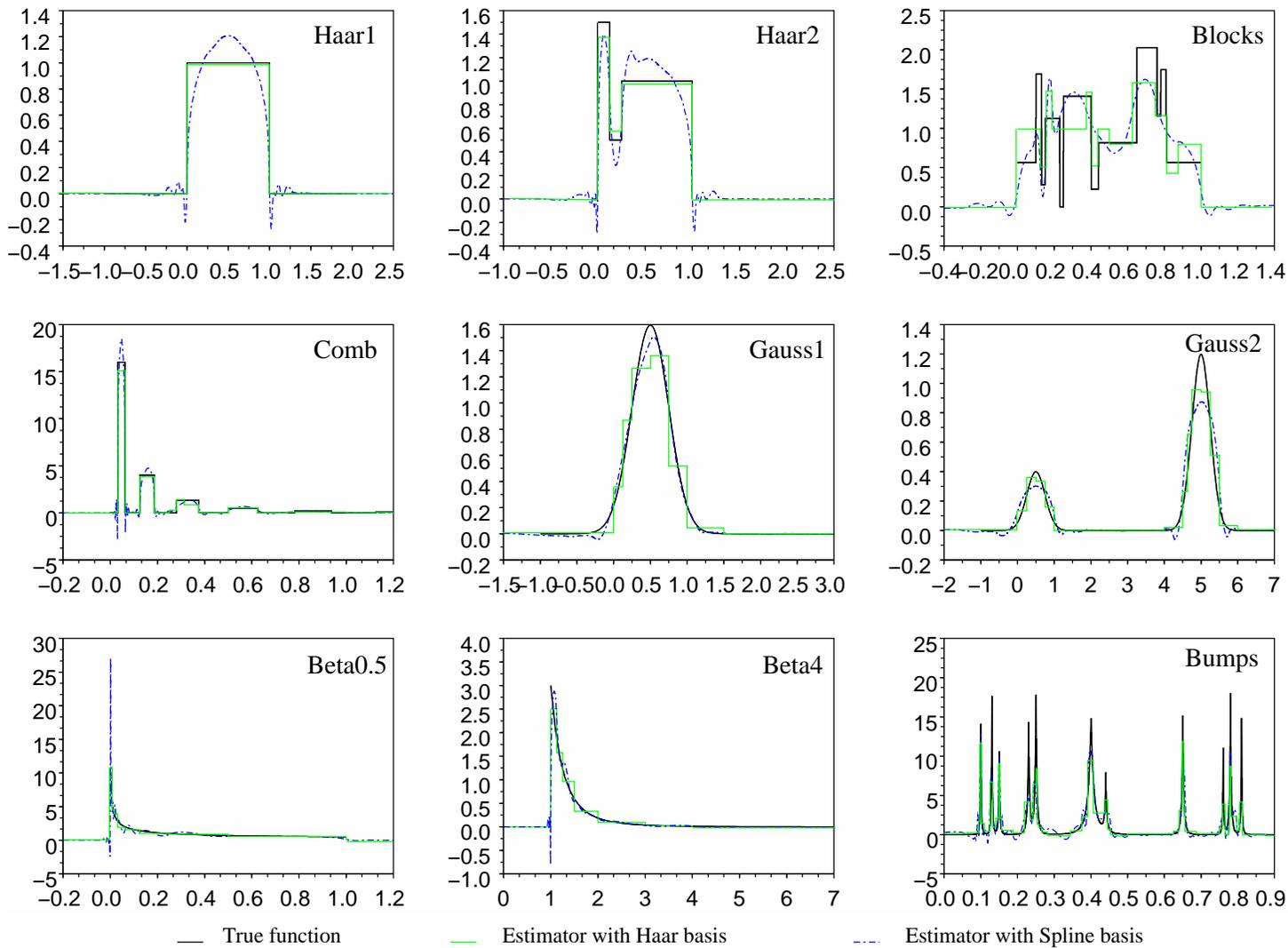}
\end{center}
\caption{Reconstructions  by  using  the  Haar   and  the  spline  bases  of  9
signals with $n=1024$, $j_0=10$ and $\gamma=1$. Top: 'Haar1', 'Haar2', 'Blocks'; Middle: 'Comb', 'Gauss1', 'Gauss2';
Bottom: 'Beta0.5', 'Beta4', 'Bumps' }\label{fig-reconstruction}
\end{figure}

The numerical performance of our procedure is first illustrated by performing it for estimating nine various signals whose definitions are given in Section \ref{defsignal}. These
functions  are  respectively   denoted  'Haar1',  'Haar2',  'Blocks',  'Comb',
'Gauss1', 'Gauss2', 'Beta0.5', 'Beta4' and 'Bumps' and
have been  chosen to represent the  wide variety of signals  arising in signal
processing.  Each  of them  satisfies  $\norm{f}_1=1$  and  can be  classified
according to the following criteria: the smoothness, the size of the support
(finite/infinite), the value of the sup norm (finite/infinite) and the shape (to be
piecewise constant  or a mixture of  peaks).  Remember that when estimating $f$, our thresholding algorithm does not use $\norm{f}_\infty$, the smoothness of $f$ and  the support of $f$ denoted $\supp(f)$ (in particular $\norm{f}_\infty$ and $\supp(f)$ can be infinite).  Simulations are performed with $n=1024$, so we observe in average $n\times\norm{f}_1=1024$ points of the underlying Poisson process. To complete the definition of  $\tilde      f_\gamma=(\tilde{f}_{n,\gamma})_n$, we rely on Theorems \ref{inegoraclelavraie}  and \ref{lower} and  we choose $j_0=\log_2(n)=10$ and $\gamma=1$ (see conclusions of Section \ref{calibrationnumerical}).  Figure \ref{fig-reconstruction}  displays intensity
reconstructions we obtain for the Haar and the spline bases.

The preliminary conclusions drawn from  Figure \ref{fig-reconstruction} are the following. As expected, a convenient choice of the wavelet system improves the reconstructions. We notice that the  estimate $\tilde  f_{n,1}$ seems to perform well for estimating the size and the location of peaks. Finally, we emphasize that the support of each signal does not play any role (compare estimation of 'Comb' which has an infinite support and  the estimation of 'Haar1' for instance). 

%%%%%%%%%%%%%%%%%%%%%%%%%%%%%%%%%%%%%%%%%%%%%%%%%
\subsection{Calibration of our procedure from the numerical point of view}\label{calibrationnumerical}
\begin{figure}[tbp]
\begin{minipage}[c]{.56\linewidth}
\hspace{-0.7cm}\includegraphics[width=\linewidth,angle=0]{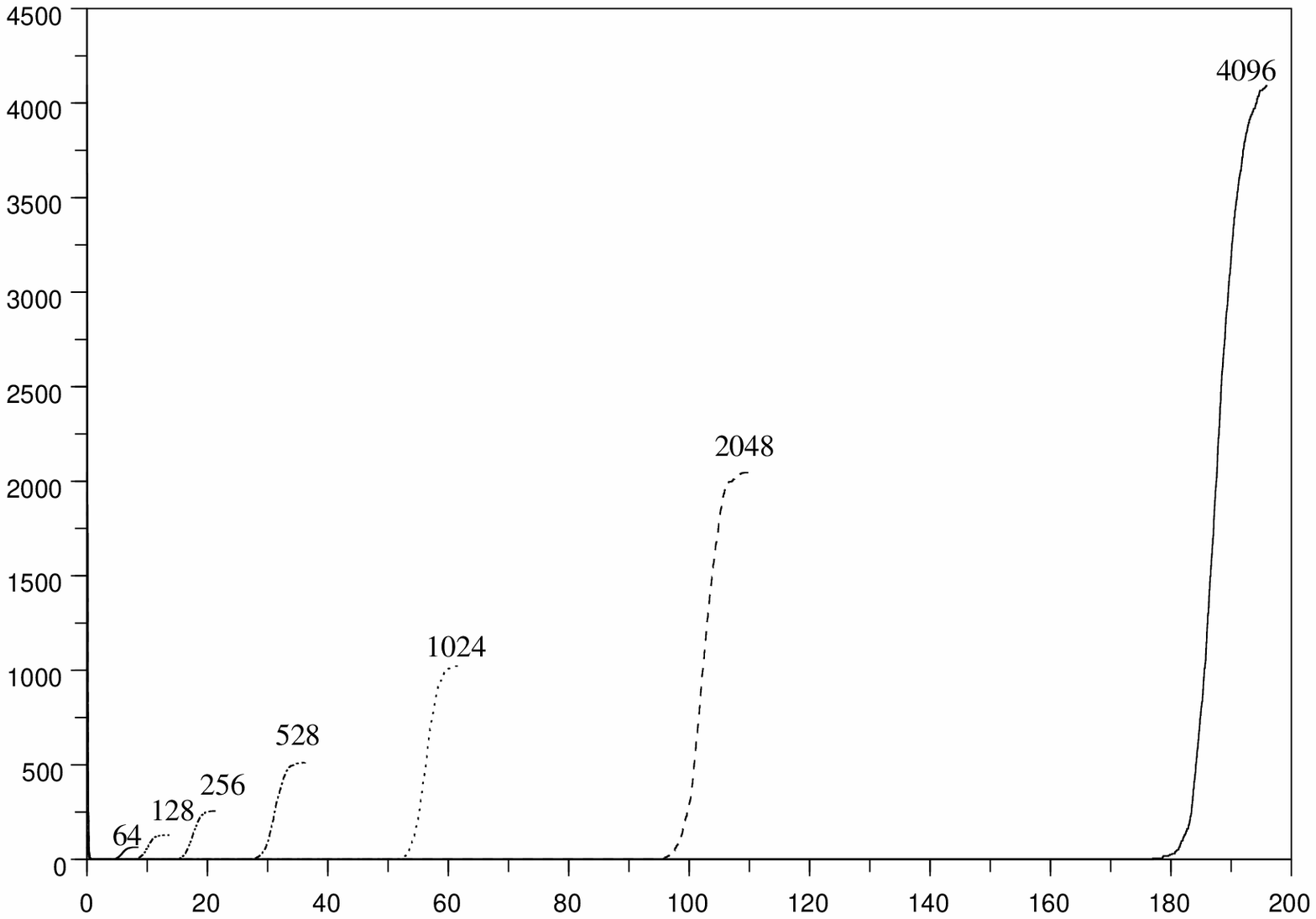}
\end{minipage} \hfill
\hspace{-1.5cm} \begin{minipage}[c]{.56\linewidth}
\includegraphics[width=\linewidth,angle=0]{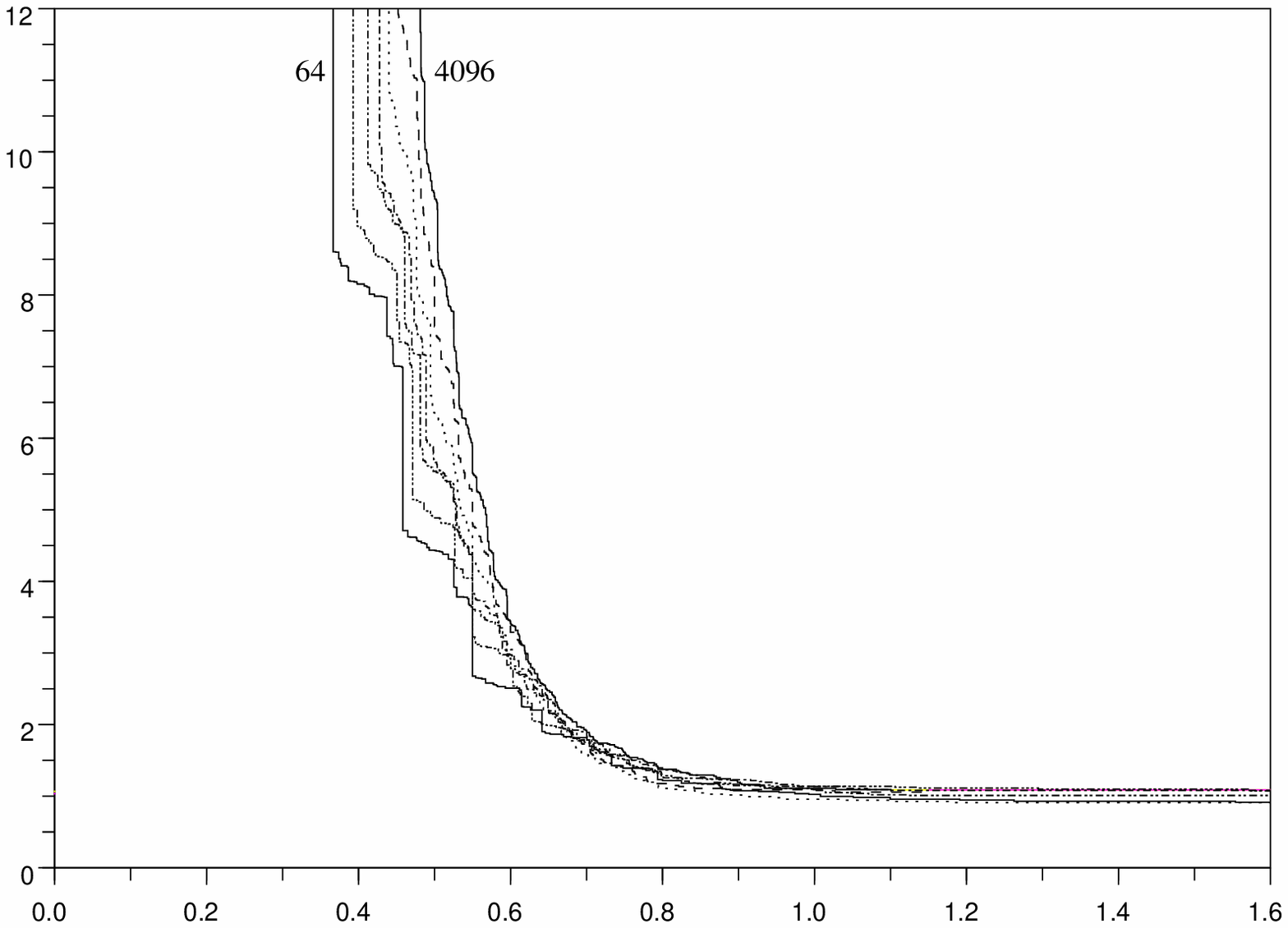}
\end{minipage}
\caption{The  function  $\gamma\to  \overline{R_n}(\gamma)$  at  two  scales  for
'Haar1'      decomposed      on      the      Haar     basis      and      for
$n\in\{64,128,256,512,1024,2048,4096\}$ with $j_0=\log_2(n)$.}\label{fig-Haar1}
\end{figure}

\begin{figure}[]
\begin{center}
\includegraphics[width=0.7\linewidth,angle=-0]{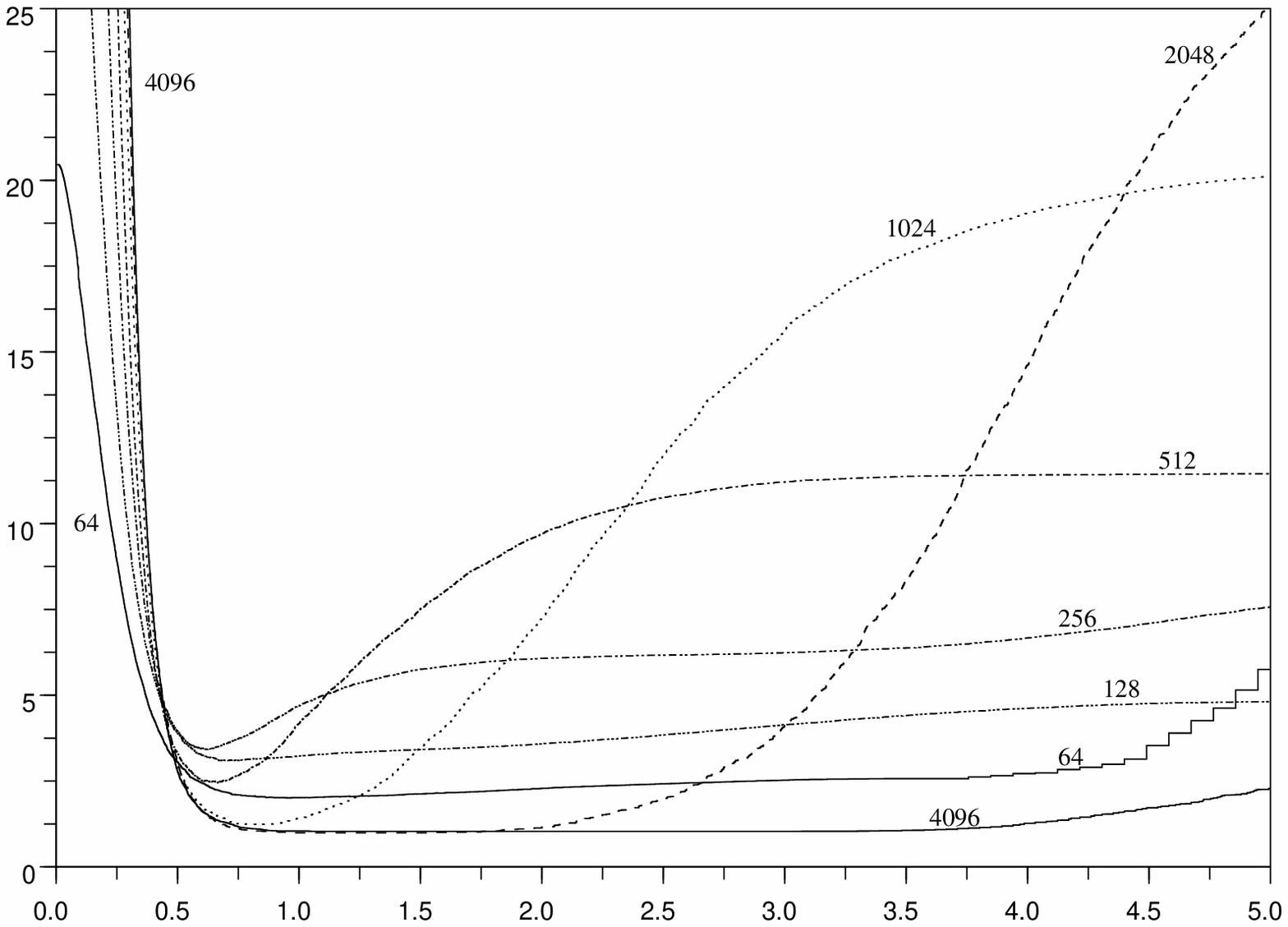}
\end{center}
\caption{The  function  $\gamma\to  \overline{R_n}(\gamma)$  for
'Gauss1'      decomposed      on      the      spline     basis      and      for
$n\in\{64,128,256,512,1024,2048,4096\}$ with $j_0=\log_2(n)$.}\label{fig-Gauss1}
\end{figure}

\begin{figure}[tbp]
\begin{minipage}[c]{.56\linewidth}
\hspace{-0.7cm}\includegraphics[width=\linewidth,angle=0]{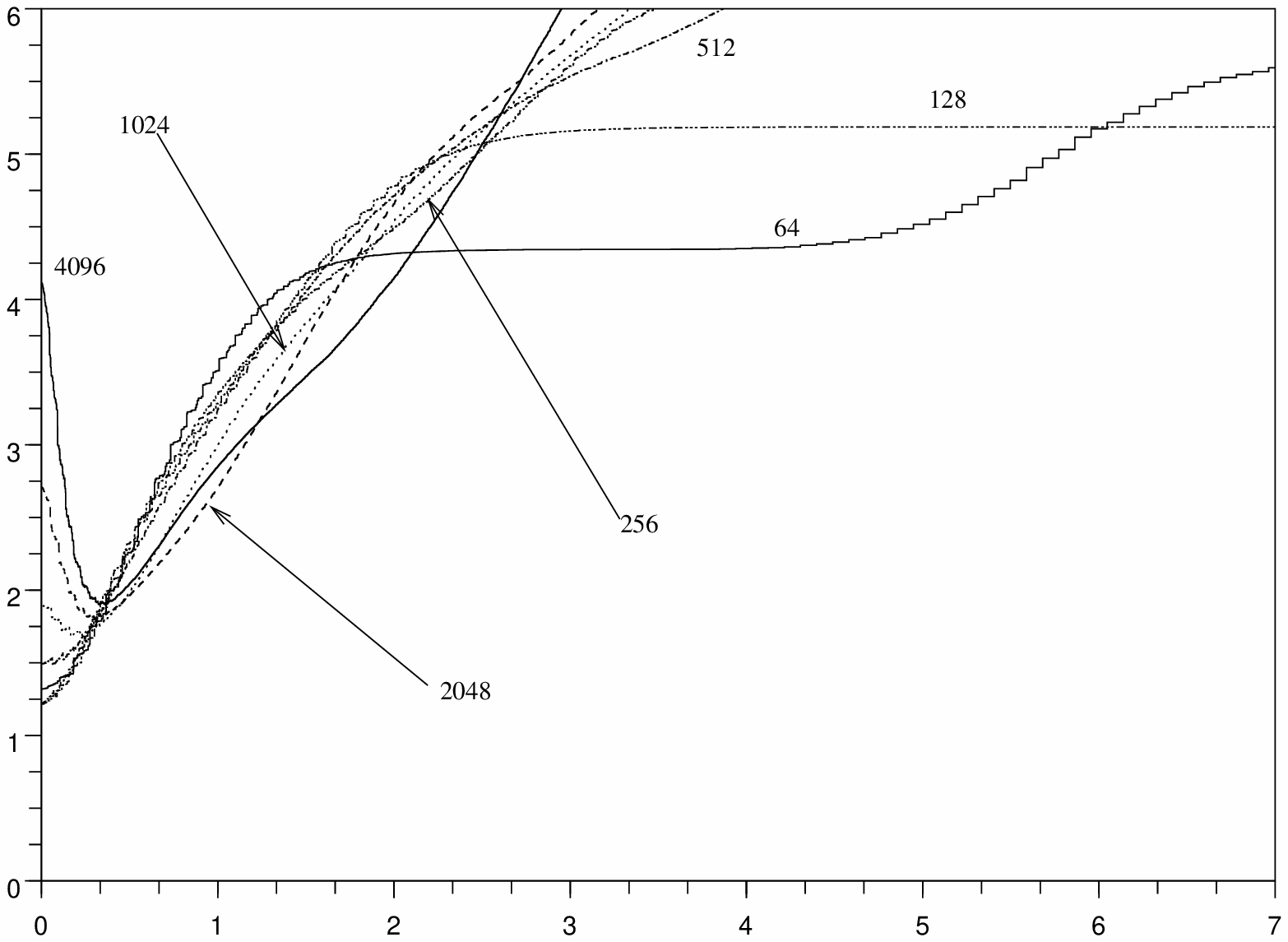}
\end{minipage} \hfill
\hspace{-1.5cm} \begin{minipage}[c]{.56\linewidth}
\includegraphics[width=\linewidth,angle=0]{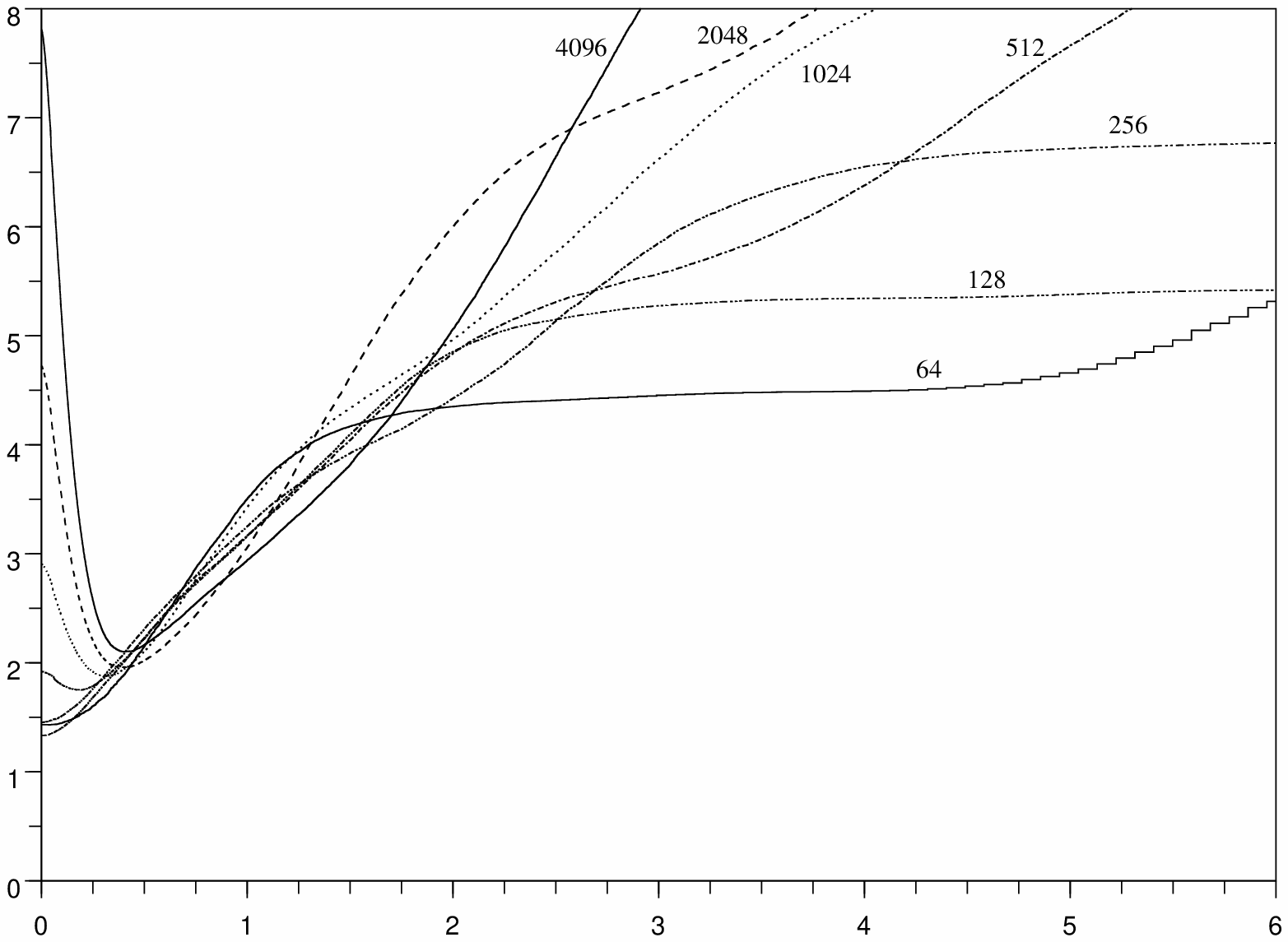}
\end{minipage}
\caption{The  function  $\gamma\to  \overline{R_n}(\gamma)$  for
'Bumps' decomposed on the Haar and the spline bases and for
$n\in\{64,128,256,512,1024,2048,4096\}$ with $j_0=\log_2(n)$.}\label{fig-Bumps}
\end{figure}

In this section, we deal with the choice of the threshold parameter
$\gamma$ in our procedures from a practical point of view. We already know that the interval $[1,12]$ is the right range for $\gamma$, theoretically speaking.  Given $n$ and a function $f$, we denote 
$R_n(\gamma)$ the ratio between the $\ell_2$-performance of our
procedure (depending on $\gamma$) and the oracle risk where the wavelet coefficients at levels $j>j_0$
are omitted. We have:
$$
R_n(\gamma)=\frac{\sum_{\la\in\Gamma_n}(\tilde\be_\la-\be_\la)^2}{\sum_{\la\in\Gamma_n}\min(\be_\la^2,V_{\la,n})}=\frac{\sum_{\la\in\Gamma_n}(\hb_\la \indic_{|\hb_\la|\geq\eta_{\la,\gamma}}-\be_\la)^2}{\sum_{\la\in\Gamma_n}\min(\be_\la^2,V_{\la,n})}. 
$$
 Of course, $R_n$ is a stepwise function and the change points of $R_n$
correspond to  the values  of $\gamma$ such  that there exists  $\lambda$ with
$\eta_{\la,\gamma}=|\hat\be_\la|$.  
The average over 1000 simulations of $R_n(\gamma)$ is computed
providing an estimation of $\E(R_n(\gamma))$. 
This average ratio, denoted $\overline{R_n}(\gamma)$ and viewed as a function of $\gamma$, is plotted for $n\in\{64,128,256,512,1024,2048,4096\}$ and for three signals considered previously:
'Haar1', 'Gauss1' and 'Bumps'. For non compactly supported signals, we need to compute an infinite number of wavelet coefficients to determine this ratio. To overcome this problem, we omit the tails of the signals and we focus our attention on an interval that contains all observations. Of course, we ensure that this approximation is negligible with respect to the values of $R_n$. As previously, we take $j_0=\log_2(n)$.
Figure  \ref{fig-Haar1}     displays
$\overline{R_n}$ for 'Haar1' decomposed on the Haar basis. The left side of Figure
\ref{fig-Haar1} gives a  general idea of the shape  of $\overline{R_n}$, while the
right side focuses on small values of $\gamma$.  
Similarly, Figures \ref{fig-Gauss1} and 
 \ref{fig-Bumps} display  $\overline{R_n}$ for 'Gauss1' decomposed on
the spline basis and for 'Bumps' decomposed on
the Haar and the spline bases. 

To discuss our results, we introduce
$$\gamma_{\min}(n)=\argmin_{\gamma>0}\overline{R_n}(\gamma).$$ 
 
For  'Haar1',
$\gamma_{\min}(n)\geq 1$ for any value of $n$ and taking
$\gamma<1$ deteriorates the performances of  the estimate.  The larger $n$, the stronger the deterioration is. Such a result was
established from the theoretical point of view in Theorem \ref{lower}. In fact,  Figure
\ref{fig-Haar1} allows to draw the following major
conclusion for 'Haar1':
\begin{equation}\label{super}
\overline{R_n}(\gamma)\approx\overline{R_n}(\gamma_{\min}(n))\approx
1
\end{equation}
 for $\gamma$ belonging to a large interval that contains the value $\gamma=1$. 
 For instance, when $n=4096$, the function $\overline{R_n}$ is close to 1 for any value of the interval $[1,177]$. So, we observe a kind of ``plateau phenomenon''. Finally, we conclude that our thresholding rule with $\gamma=1$ performs very well since it achieves the same performance as the oracle estimator.

For 
'Gauss1', $\gamma_{\min}(n)\geq 0.5$ for any value of $n$. Moreover, as soon
as $n$ is large enough, the oracle ratio for $\gamma_{\min}(n)$ is close to $1$. Besides, when $n\geq 2048$,
as for  'Haar1', $\gamma_{\min}(n)$  is larger than  $1$. We observe the
``plateau phenomenon'' as well and as for 'Haar1', the size of the plateau increases when $n$ increases.  This can be
explained by the following important  property of 'Gauss1'\pa{:} 'Gauss1' can be well approximated by a
finite combination of the atoms of the spline basis. So, we have the strong
impression  that  the  asymptotic  result  of  Theorem  \ref{lower}  could  be
generalized for  the spline  basis.  

Conclusions for 'Bumps' are very different. Remark that this irregular signal
has many significant wavelet coefficients at
high resolution levels whatever  the basis. We have $\gamma_{\min}(n)<0.5$ for
each  value of  $n$. Besides,  $\gamma_{\min}(n)\approx 0$  when  $n\leq 256$,
which means that all the coefficients until  $j=j_0$ have to be kept to obtain the
best estimate. So, the parameter $j_0$ plays an essential role and has to be
well  calibrated   to  ensure  that   there  are  no   non-negligible  wavelet
coefficients for  $j>j_0$.  Other differences  between Figure \ref{fig-Haar1}
(or Figure \ref{fig-Gauss1}) and Figure \ref{fig-Bumps} have to be emphasized.  For
'Bumps', when $n\geq 512$, the minimum of $\overline{R_n}$ is well localized,
there is no plateau anymore and $\overline{R_n}(1)>2$. Note that $\overline{R_n}(\gamma_{\min}(n))$ is larger than 1. 

\begin{figure}[tbp]
\begin{center}
\includegraphics[width=0.7 \linewidth,angle=-0]{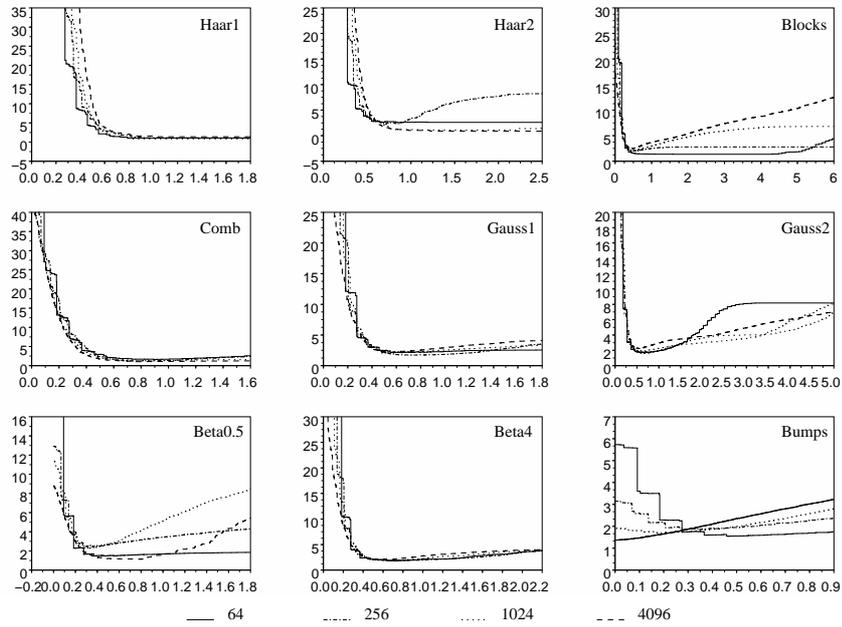}
\end{center}
\caption{Average over  100 iterations of the  function $R_n$
for signals decomposed on the Haar basis and for
$n\in\{64,256,1024,4096\}$ with $j_0=10$.}\label{tous-1}
\end{figure}

\begin{figure}[htbp]
\begin{center}
\includegraphics[width=0.7\linewidth,angle=-0]{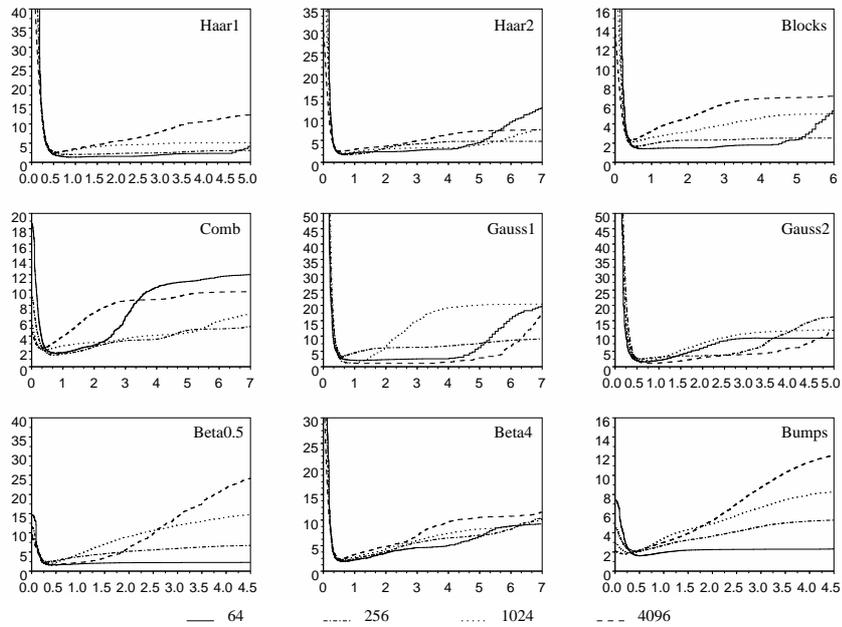}
\end{center}
\caption{Average over  100 iterations of the  function $R_n$
for signals decomposed on the spline basis and for
$n\in\{64,256,1024,4096\}$ with $j_0=10$.}\label{tous-2}
\end{figure}

Previous preliminary conclusions show that  the ideal  choice  for  $\gamma$  and  the
performance of  the thresholding rule  highly depend on the  decomposition of
the signal on the wavelet basis. Hence, in the sequel, we have decided to take $j_0=10$ 
 for any value of $n$ so
that the decomposition on the basis is not too coarse. To extend previous results, Figures
\ref{tous-1} and  \ref{tous-2} display the  average of the function  $R_n$ for
the signals 'Haar1',  'Haar2',  'Blocks',  'Comb',
'Gauss1', 'Gauss2', 'Beta0.5', 'Beta4' and 'Bumps' with $j_0=10$. For the sake of brevity, we only consider the values
$n\in\{64,256,1024,4096\}$ and the average of
$R_n$ is performed over 100 simulations.
 Figure \ref{tous-1} gives the results
obtained for the  Haar basis and Figure \ref{tous-2} for  the spline basis.
This study allows to draw conclusions with respect to the issue of calibrating $\gamma$ from the numerical point of view. To present them, let us introduce two classes of functions. 

The first class is the class of
signals that only have negligible coefficients at high levels of resolution. The wavelet basis is well adapted to the signals of this class that
contains 'Haar1', 'Haar2' and 'Comb' for the Haar basis and 'Gauss1' and 'Gauss2' for the spline
basis.  For such  signals, the estimation problem  is close  to  a parametric
problem. In  this  case,  the  performance  of the  oracle  estimate can  be
achieved at least for $n$ large enough and (\ref{super}) is true for $\gamma$ belonging to a large interval that contains the value $\gamma=1$. These numerical conclusions strengthen and generalize theoretical conclusions of Section \ref{penaltyterm}.

The second class of functions is the class of irregular signals with significant wavelet coefficients at high
resolution levels. For such signals $\gamma_{\min}(n)<0.8$ and
there is no ``plateau'' phenomenon (in particular, we do not have $\overline{R_n}(1)\simeq\overline{R_n}(\gamma_{\min}(n))$). 

Of course, estimation is easier and performances of our procedure are better when the signal belongs to the first class. But in practice, it is hard to  choose a wavelet system such that the intensity to be estimated satisfies this property. However, our study allows to use the following simple rule. If the practitioner has no idea of the ideal wavelet basis to use, he should perform the thresholding rule with $\gamma=1$ (or $\gamma$ slightly larger than 1) that leads to convenient results whatever the class the signal belongs to.
%%%%%%%%%%%%%%%%%%%%%%%%%%%%%%%%%
\subsection{Comparisons with classical procedures}\label{comparaison}
\begin{figure}[tpb]
\begin{center}
\includegraphics[width=0.6\linewidth,angle=0]{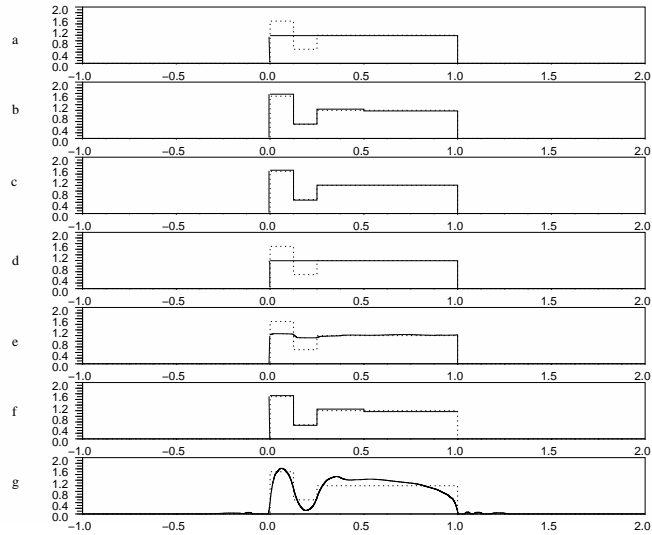}
\end{center}
\caption{Reconstructions of 'Haar2' with $n=1024$. (a) ANSCOMBE-UNI; (b) CORRECTED; (c) RAND-THRESH-HAAR; (d) ANSCOMBE-UNI-TI; (e)  CORRECTED-TI; 
(f) FREE-DEGREE; (g) RAND-THRESH-SPLINE.}\label{fig-reconstruction-Haar2}
\end{figure}

\begin{figure}[tpb]
\begin{center}
\begin{minipage}[c]{.5\linewidth}
\hspace{.5cm}\includegraphics[width=\linewidth,angle=0]{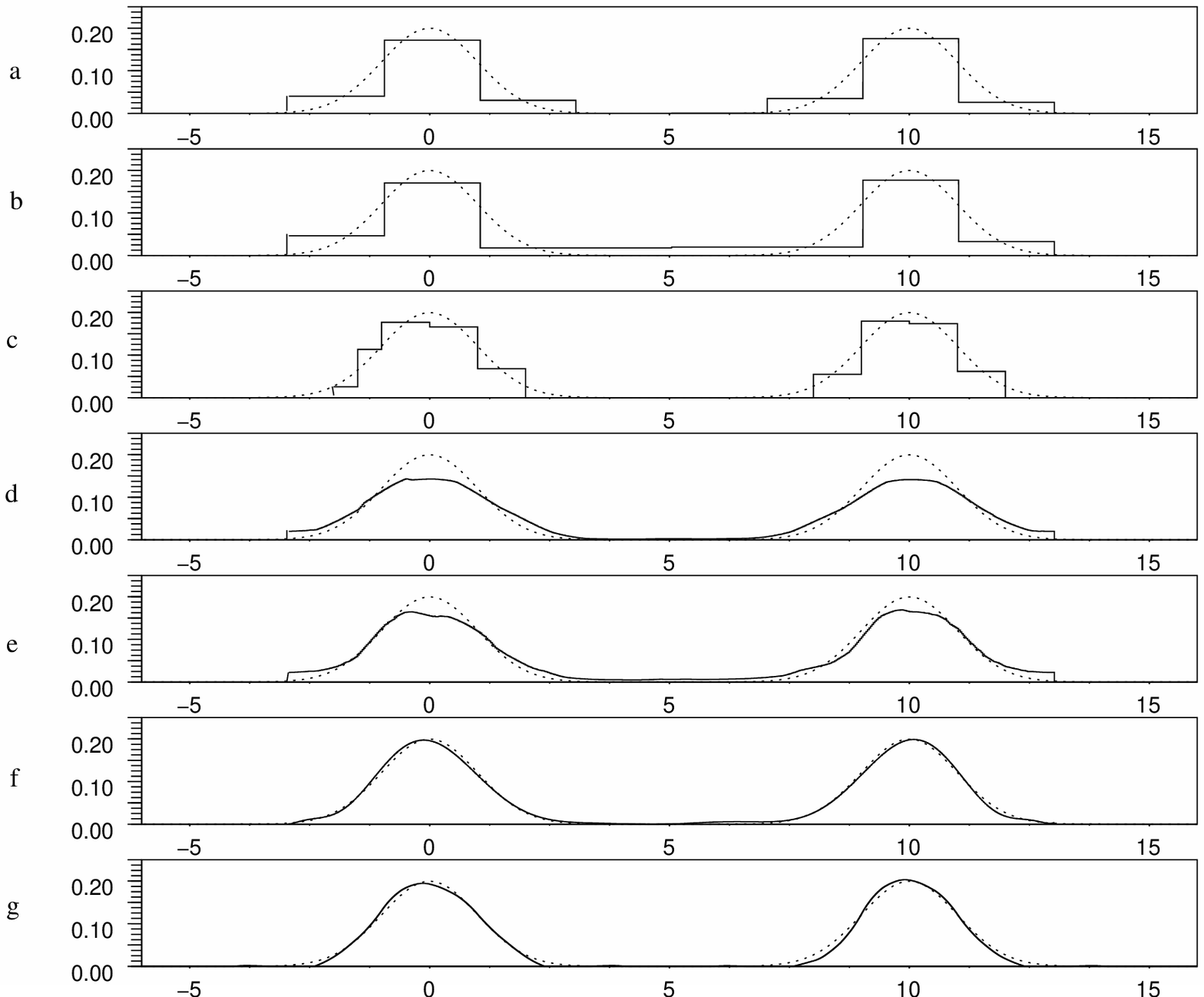}
\end{minipage} \hfill
\hspace{-3.5cm} \begin{minipage}[c]{.5\linewidth}
\includegraphics[width=\linewidth,angle=0]{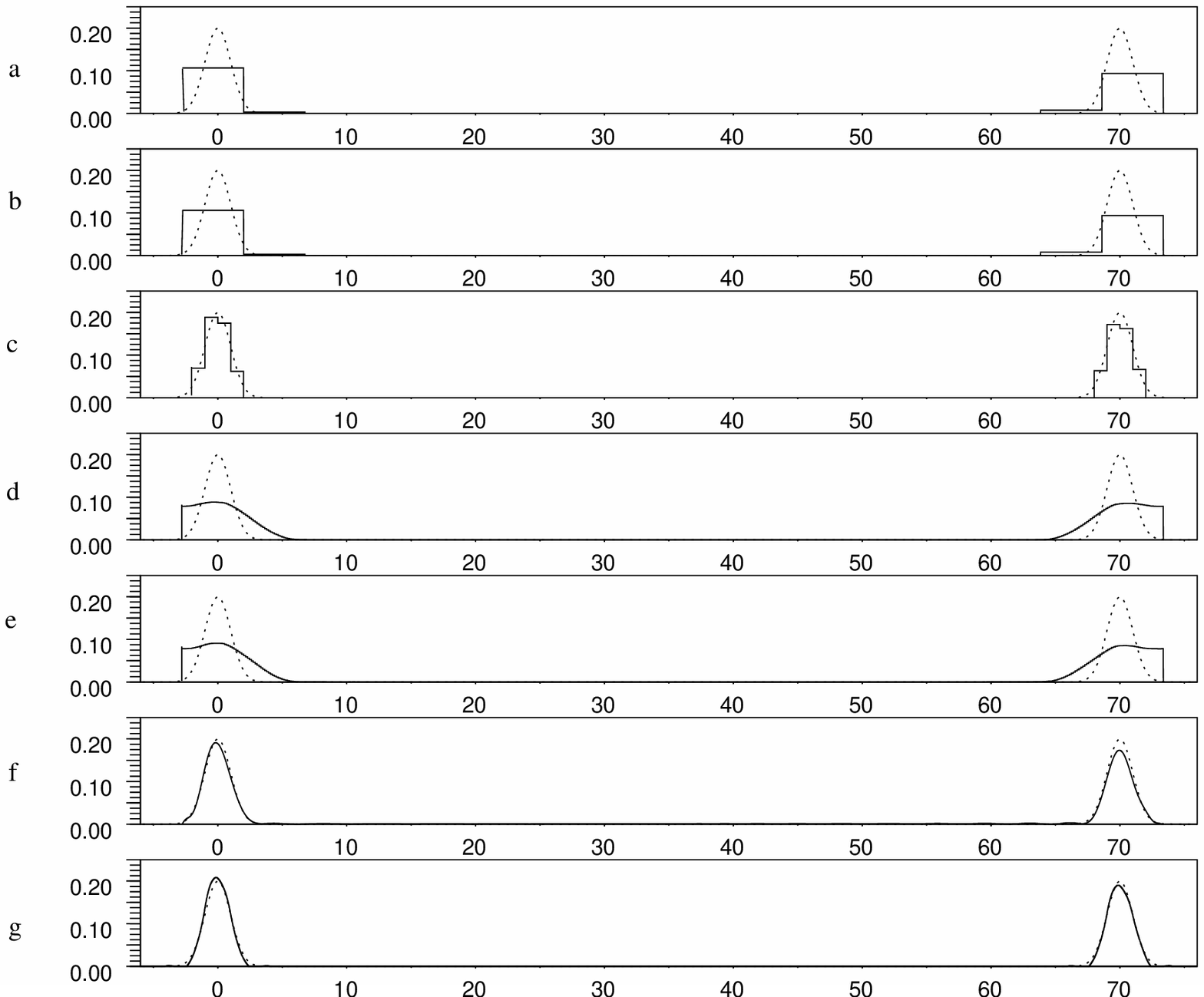}
\end{minipage}
\end{center}
\caption{Reconstructions of $f_d$ with $n=1024$ (left: $d=10$, right $d=70$). (a) ANSCOMBE-UNI; (b) CORRECTED; (c) RAND-THRESH-HAAR; (d) ANSCOMBE-UNI-TI; (e) CORRECTED-TI; 
(f) FREE-DEGREE; (g) RAND-THRESH-SPLINE.}\label{fig-reconstruction-Gauss-ecart}
\end{figure}

\begin{figure}[tpb]
\begin{center}
\includegraphics[width=0.6\linewidth,angle=-0]{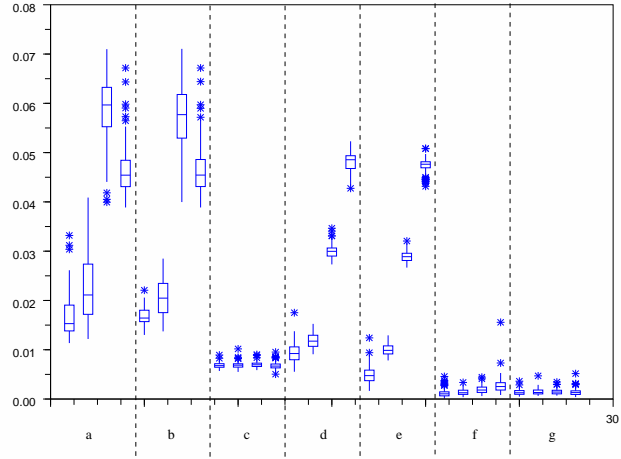}
\end{center}
\caption{Mean square error over 100 simulations of the different methods with $n=1024$. From left to right: 10, 30, 50 and 70. (a): ANSCOMBE-UNI; (b): CORRECTED
; (c): RAND-THRESH-HAAR; (d): ANSCOMBE-UNI-TI; (e) : CORRECTED-TI; 
(f): FREE-DEGREE; (g): RAND-THRESH-SPLINE.}\label{fig-moustache}
\end{figure}

Now, let us compare our procedure with classical ones. We first consider the methodology based on the Anscombe transformation of Poisson type observations  (see \cite{ans}). This preproprecessing yields Gaussian data with a constant noise level close to 1. Then, universal wavelet thresholding proposed by Donoho and Johnstone \cite{dojo} is applied with the Haar basis. Kolaczyk corrected  this standard algorithm  for burst-like Poisson data. He proposed to use Haar wavelet thresholding directly on the binned data with
especially calibrated thresholds (see \cite{kolastro} and \cite{kol}). In the sequel, these algorithms are respectively denoted ANSCOMBE-UNI and CORRECTED. We briefly mention that CORRECTED requires the knowledge of a so-called background rate that is empirically estimated in our paper (note however that CORRECTED heavily depends on the precise knowledge of the background rate as shown by the extensive study of Besbeas, de Feis and Sapatinas \cite{bfs}). One can combine the wavelet transform and translation invariance to eliminate the shift dependence of the Haar basis.  When ANSCOMBE-UNI and CORRECTED are combined with translation invariance, they are respectively denoted ANSCOMBE-UNI-TI and CORRECTED-TI in the sequel. Finally, we consider the penalized piecewise-polynomial rule proposed by Willett and Nowak \cite{wn} (denoted FREE-DEGREE in the sequel) for multiscale Poisson intensity estimation. Unlike our estimator, the knowledge of the support of $
 f$ is essential to perform all these procedures that will be sometimes called ``support-dependent strategies'' along this section. We first consider estimation of the signal 'Haar2' supported by $[0,1]$ for which reconstructions  with $n=1024$ are proposed in Figure \ref{fig-reconstruction-Haar2} where we have taken the positive part of each estimate. For ANSCOMBE-UNI, CORRECTED and their counterparts based on translation invariance, the finest resolution level for thresholding is chosen to give good overall performances. For our random thresholding procedures, respectively based on the Haar and spline bases and respectively denoted RAND-THRESH-HAAR and RAND-THRESH-SPLINE, we still use $\gamma=1$ and $j_0=\log_2(n)=10$.  We note that for the setting of  Figure \ref{fig-reconstruction-Haar2}, translation invariance oversmooths estimators. Furthermore, comparing (a), (b) and (c), we observe that universal thresholding is too conservative. Our procedure works well provided the 
 Haar basis is chosen, whereas FREE-DEGREE automatically selects a piecewise constant estimator.
Now, let us consider a non-compactly supported signal based on a mixture of two Gaussian densities. We denote $d$ the distance between modes of these Gaussian densities, so the intensity associated with this signal is
$$f_d(x)=\frac{1}{2}\left(\frac{1}{\sqrt{2\pi}}\exp\left(-\frac{x^2}{2}\right)+\frac{1}{\sqrt{2\pi}}\exp\left(-\frac{(x-d)^2}{2}\right)\right)$$
and we take $n=1024$. To apply support-dependent strategies, we consider the interval given by the smallest and the largest observations and data are  first rescaled to be supported by the interval $[0,1]$. Reconstructions with $d=10$ and $d=70$ are given in Figure \ref{fig-reconstruction-Gauss-ecart}.
RAND-THRESH-HAAR outperforms ANSCOMBE-UNI and CORRECTED but all these procedures are too rough. To some extent, it is also true for ANSCOMBE-UNI-TI and CORRECTED-TI even if translation invariance improves the corresponding reconstructions. This is not the case for RAND-THRESH-SPLINE and  FREE-DEGREE. When $d=70$, performances of  all the support-dependent strategies
deteriorate, which illustrates the harmful role of the support. 
 In particular, procedures based on the translation invariance principle which periodizes the data, deal with the two main parts of the signal as if they were close to each other, \pa{they} are \pa{consequently} quite inadequate. The worse performances of FREE-DEGREE for $d=70$ could be expected since its theoretical performances are established under the strong assumption that the signal is bounded from below on its \pa{(known)} support. To strengthen these results and to show the influence of the support, we compute the  mean square error over
100 simulations for each method and we provide
the corresponding  boxplots given in Figure \ref{fig-moustache} associated with $f_d$ when $d\in\{10,30,50,70\}.$ 
Note that when $d$ increases, unlike the other algorithms, performances of our thresholding rule based either on the Haar or on the spline basis are remarkably stable. In particular,  for $d=70$, RAND-THRESH-SPLINE outperforms \pa{all} the other algorithms.
\pa{Note also the very bad performances of ANSCOMBE-UNI and CORRECTED for $d=50$ due to the inadequacy between the way the data are binned and the distance $d$.} 

The main conclusions of this short study are the following. We note that the estimate proposed in this paper outperforms ANSCOMBE-UNI 
and CORRECTED (compare (a), (b) and (c)), showing that the data-driven calibrated threshold proposed in (\ref{defthresh}) improves classical ones. In particular, classical methods highly depend on the way data are binned and on the choice of resolutions levels where coefficients are thresholded, whereas our methodology only depends on $\gamma$ and on $j_0$ for which we propose to take systematically $\gamma=1$ and $j_0=\log_2 (n)$.  However, unlike FREE-DEGREE, we have to choose a convenient wavelet basis for decomposing \pa{the} signals. Finally, the support, if too large, can play a harmful role whenever the method needs to rescale the data. This is not the case for the method presented in this paper, which explains the robustness of our procedures with respect to the support issue.
%%%%%%%%%%%%%%%%%%%%%%%%%%%%%%%%%%%%%%%%%%%
%%%%%%%%%%%%%%%%%%%%%%%%%%%%%%%%%%%%%%%%%%%
\section{Proofs of the results}\label{proofs}
%Using (\ref{equnormes}), without loss of generality,  the results are  established by  using the  $\ell_2$-norm of  coefficients instead  of the functional $\L_2$-loss.
%%%%%%%%%%%%%%%%%%%%%%%%%%%%%%%%%%%%%%%%%%%
\subsection{Proof of Proposition \ref{Propsharp}}\label{preuvepropsharp}
The first point is obvious. For the second point, first, let us take $f\in{\cal F}$. We can write $f=\sum_{\la\in
\Lambda_1}\be_\la\tilde\p_\la$,
where  $$\Lambda_1=\{\la:\ \be_\la\not=0\}$$  is finite. Since $\be_\la\not=0$
implies  $F_\la>0$, we have $$\min_{\la\in\Lambda_1}F_\la>0.$$  So,  $f$ belongs  to
${\cal F}_{n}(R)$ for $n$ and $R$ large enough.\\ Conversely,  if       $f=\sum_{\la\in
\Lambda}\be_\la\tilde\p_\la$ belongs  to ${\mathcal F}_n(R)$ for some  $n$ and some $R>0$ and if
$f$ has an infinite number of non-zero wavelet coefficients, then there is an
infinite  number of  indices $\lambda=(j,k)$  such that  $$F_\la=F_{j,k}\geq\frac{(\ln n)(\ln\ln
n)}{n}.$$   So,   either   for   any  arbitrary large  $j$,  there   exists   $k$   such   that 
$$\frac{(\ln n)(\ln\ln n)}{n}\leq F_{j,k}\leq \norm{f}_\infty |\supp(\p_{j,k})|=\norm{f}_\infty 2^{-j},$$
so $f\not\in \L_\infty(R)$ or there exists
$j$  such that $\sum_kF_{j,k}=+\infty$  and $f\not\in  \L_1(R)$ (see (\ref{Fla})). This  cannot occur
since $f\in {\mathcal F}_n(R)$. This \pa{concludes the proof of Proposition 1}.
%%%%%%%%%%%%%%%%%%%%%%%%%%%%%%%%%%%%%%%%%%%
\subsection{Proof of Theorem \ref{classFn}}\label{sec:proofclassFn}
We first state the following lemma established in \cite{Poisson_minimax} where it is used to derive Theorem \ref{inegoraclelavraie}.  For the sake of exhaustiveness, the proof of Lemma \ref{lemmesharp} is recalled in section \ref{preuvelemmesharp}.
\pa{\begin{Lemma}\label{lemmesharp}
For all $\kappa$ such that $\gamma^{-\frac{1}{2}}<\kappa<1$, there exists a positive constant $K$ depending on $\gamma,~\kappa$ and $\norm{f}_1$ such that
\begin{equation}
\E\normp{\tilde{f}^H_{n,\gamma}-f}^2\leq \left(\frac{1+\kappa^2}{1-\kappa^2}\right)\inf_{m\subset
  \Ga_n}\left\{\frac{1+\kappa^2}{1-\kappa^2}\sum_{\la\not \in
  m}\be_\la^2+\frac{1-\kappa^2}{\kappa^2}\sum_{\la\in
  m}\E(\hb_\la-\be_\la)^2+\sum_{\la \in m}\E(\eta_{\la,\gamma}^2)\right\}+\frac{K}{n},\nonumber
\end{equation}
where we denote by $m$ any possible subset of indices $\lambda$.
\end{Lemma}}
First, we give an upper bound for $\E(\eta_{\la,\gamma}^2)$.
For any $\delta>0$, 
$$
\E(\eta_{\la,\gamma}^2)\leq (1+\delta) 2 \ga \ln n \E(\tilde{V}_{\la,n})+ (1+\delta^{-1})\left(\frac{\ga \ln n}{3n}\right)^2\norm{\p_\la}_\infty^2.
$$
Moreover, $$\E(\tilde{V}_{\la,n})\leq (1+\delta) V_{\la,n} +(1+\delta^{-1}) 3
\gamma \ln n \frac{\norm{\p_\la}_\infty^2}{n^2}.$$
So,
\begin{equation}\label{majothresh}
\E(\eta_{\la,\gamma}^2)\leq (1+\delta)^2 2 \ga \ln n  V_{\la,n} + \Delta(\delta) \left(\frac{\ga \ln n}{n}\right)^2\norm{\p_\la}_\infty^2,
\end{equation}
with $\Delta(\delta)$ a constant depending only on $\delta$.
Now, let us choose the parameter $\gamma$ in an optimal way. The main terms in the upper bound given by the lemma are the first and third ones. So,  we choose $\kappa^2$ close to $\gamma^{-1}$ as required by the assumptions to the lemma and we fix $\gamma$ such that
$$\left(\frac{1+\kappa^2}{1-\kappa^2}\right)^2\approx\left(\frac{\gamma+1}{\gamma-1}\right)^2\quad\mbox{and}\quad 2\gamma\left(\frac{1+\kappa^2}{1-\kappa^2}\right)\approx \frac{2(\gamma^2+\gamma)}{\gamma-1}$$ 
are as small as possible. %Even if this disadvantages the term $\sum_{\la\not \in  m}\be_\la^2$, 
We first minimize $\frac{2(\gamma^2+\gamma)}{\gamma-1}$ so we choose $\ga=1+\sqrt{2}$. Now, we set $\kappa=\sqrt{0.42}\approx (1+\sqrt{2})^{-1/2}$. Then, with  $\delta>0$ such that $$(1+\delta)^2=11.822(1-\kappa^2)(2\ga(1+\kappa^2))^{-1}\simeq 1.00006,$$ we obtain
\begin{eqnarray*}
\E\normp{\tilde f_{n,\ga}^H-f}^2&\leq& \inf_{m\subset \Ga_n} \left\{6\sum_{\la\not\in m} \be_\la^2+\sum_{\la\in m} (3.4 + 11.822\ln n) V_{\la,n} + \Delta'\sum_{\la\in m}\left(\frac{\ln n \norm{\p_\la}_\infty}{n}\right)^2\right\}+ \frac{K}{n},
\end{eqnarray*}   
where $$\Delta'=\Delta(\delta)\gamma^2(1+\kappa^2)(1-\kappa^2)^{-1}.$$ 
Let $n$ and $R>0$ be fixed and let $f\in {\cal F}_n(R)$.
Assume that $\be_\la\not=0$. 
%Since $$|\be_\la|\leq 2^{\max\left(0,\frac{j}{2}\right)}F_\la,$$where $\la=(j,k)$, then $F_\la\not=0$. 
In this case, 
$$F_\la\geq \frac{(\ln n)(\ln\ln n)}{n}.$$ %and using
\pa{But} 
$$F_\la  \leq 2^{-\max(j,0)}\norm{f}_{\infty}\leq 2^{-\max(j,0)}R$$ for $\la=(j,k)$. \pa {So} % we have
$2^j\leq 2^{j_0}$ \pa{holds} for $n$ large enough and $\la$ belongs to $\Ga_n$. Finally, we conclude that  $\be_\la\not=0$ implies   $\la    \in\Ga_n$. Now,   take   $$m=\{\la\in   \Ga_n:\quad
\be_\la^2>V_{\la,n}\}.$$ If $m$ is empty, then $\be_\la^2=\min(\be_\la^2,V_{\la,n})$ for every
$\la\in\Ga_n$. Hence
$$\E\normp{\tilde f_{n,\ga}^H-f}^2\leq 6\sum_{\la\in \Ga_n}
\min(\be_\la^2,V_{\la,n}) + \frac{K}{n}$$
and Theorem \ref{classFn} is proved. If $m$ is not empty, with $\la=(j,k)$, 
$$V_{\la,n}=\frac{2^{\max(j,0)} F_\la }{n}=	 \frac{\norm{\p_\la}_\infty^2 F_\la}{n}.$$
 Hence,  for  all  $n$,  if  $\la\in  m$, then $\be_\la\not=0$ and  $$V_{\la,n}  \ln  n  \geq  \frac{(\ln  n
)^{2}(\ln\ln n) \norm{\p_\la}_\infty^2}{n^2} $$
and if $n$ is large enough, $$0.1\,\ln n \sum_{\la \in m} V_{\la,n} \geq  \Delta'\sum_{\la \in m} \left(\frac{\ln n \norm{\p_\la}_\infty}{n}\right)^2 +3.4\sum_{\la \in m} V_{\la,n}.$$ Theorem \ref{classFn} is proved since for $n$ large enough (that depends on $R$), we obtain:
$$\E\normp{\tilde f_{n,\ga}^H-f}^2\leq 6\sum_{\la\not\in m} \be_\la^2+11.922\,\ln n\sum_{\la\in m} V_{\la,n}+\frac{K}{n}\leq 12\,\ln n\left(\sum_{\la\not\in m} \be_\la^2+\sum_{\la\in m} V_{\la,n}+\frac{1}{n}\right).$$

%%%%%%%%%%%%%%%%%%%%%%%%%%%%%%%%%%%%%%%%%%%%%%%
\subsection{Proof of Theorem \ref{lower}}\label{prooflower}
Let $\ga<1$. Note that for all $\e>0$,
\begin{equation}\label{inegseuil}
\sqrt{2\gamma \hat{V}_{\la,n}\ln n }+\frac{\gamma\ln
    n}{3n}\norm{\p_\la}_\infty\leq\eta_{\la,\ga}\leq \eta'_{\la,\gamma}:=\sqrt{2\ga(1+\e)\ln(n)\hat{V}_{\la,n}}+ \frac{\gamma\ln(n)\norm{\p_\la}_{\infty}}{n}w_\e,
\end{equation}
where $w_\e=\sqrt{\e^{-1}+6}+1/3$ depends only on $\e$. We choose $\e$ such that $\ga'=\ga(1+\e)<1$.
Let $\al>1$ and $n$ be fixed. We set $j$ the positive integer such that
$$\frac{n}{(\ln n)^\al}\leq 2^j <\frac{2n}{(\ln n)^\al}.$$ 
For all $k\in\{0,...,2^j-1\}$, we define $$N^+_\jk = \int_{k2^{-j}}^{(k+\frac{1}{2})2^{-j}}dN\quad\mbox{ and }\quad N^-_\jk =
\int_{(k+\frac{1}{2})2^{-j}}^{(k+1)2^{-j}}dN.$$ These variables are i.i.d. random Poisson
variables of parameter $\mu_{n,j} = n2^{-j-1}$.
Moreover,
$$\hb_\jk = \frac{2^{\frac{j}{2}}}{n}(N^+_\jk-N^-_\jk) \quad\mbox{ and }\quad
\hat{V}_{(\jk),n} = \frac{2^{j}}{n^2}(N^+_\jk+N^-_\jk).$$
Hence,
\begin{eqnarray*}
\E(\normp{\tilde{f}_{n,\gamma}^H-f}^2)&\geq & \sum_{k=0}^{2^j-1} \E\left(\hb_{j,k}^2\indic_{|\hb_{j,k}|>\eta_{\la,\ga}}\right)\\
&\geq & \sum_{k=0}^{2^j-1} \E\left(\hb_{j,k}^2\indic_{|\hb_{j,k}|>\eta'_{\la,\gamma}}\right)\\
&\geq & \sum_{k=0}^{2^j-1}
\frac{2^j}{n^2}\E\left((N^+_\jk-N^-_\jk)^2 \indic_{|N^+_\jk-N^-_\jk|\geq
\sqrt{2\gamma'\ln (n) (N^+_\jk+N^-_\jk)}+\ln(n)\gamma w_\e}\right).
\end{eqnarray*}
Let $u_n$ be a bounded sequence that will be fixed later such that $u_n\geq\gamma w_\e$. We set $$v_{n,j}= \left(\sqrt{4\gamma'\ln(n) \tilde\mu_{n,j}}+\ln(n) u_n\right)^2$$
where $\tilde\mu_{n,j}$ is the largest integer smaller that $\mu_{n,j}$.
Note  that if $$N^+_\jk=\tilde\mu_{n,j} + \frac{\sqrt{v_{n,j}}}{2}\quad\mbox{ and }\quad N^-_\jk= \tilde\mu_{n,j} -\frac{\sqrt{v_{n,j}}}{2},$$ then $$|N^+_\jk-N^-_\jk|= \sqrt{2\gamma'\ln(n) (N^+_\jk+N^-_\jk)}+\ln(n) u_n.$$ 
Let $N^+$ and $N^-$ be two independent Poisson variables of parameter
$\mu_{n,j}$. Then,
$$\E(\normp{\tilde{f}_{n,\gamma}^H-f}^2)\geq 
\frac{2^{2j}}{n^2}v_{n,j}
\P\left(N^+=\tilde\mu_{n,j} +\frac{\sqrt{v_{n,j}}}{2}\mbox{ and }N^-=\tilde\mu_{n,j} - \frac{\sqrt{v_{n,j}}}{2}\right).
$$
Note that $$\frac{1}{4}(\ln n)^\al-1<\tilde\mu_{n,j}\leq \mu_{n,j}\leq \frac{1}{2}(\ln n)^\al$$ and $$\lim_{n\to +\infty}\frac{\sqrt{v_{n,j}}}{\mu_{n,j}}=\lim_{n\to +\infty}\frac{\sqrt{v_{n,j}}}{\tilde\mu_{n,j}}=0.$$  So, we set $$l_{n,j}=\tilde\mu_{n,j}+\frac{\sqrt{v_{n,j}}}{2}\quad\mbox{ and }\quad m_{n,j}=\tilde\mu_{n,j}-\frac{\sqrt{v_{n,j}}}{2}$$ that go to
$+\infty$ with $n$. Now, we take a bounded sequence $u_n$  such that for any $n$, $\frac{\sqrt{v_{n,j}}}{2}$ is an integer and $u_n\geq\gamma w_\e$. 
Hence by the Stirling formula, 
\begin{eqnarray*}
\E(\normp{\tilde{f}_{n,\gamma}^H-f}^2)&\geq& \frac{v_{n,j}}{(\ln n)^{2\al}}
\P\left(N^+=\tilde\mu_{n,j}+\frac{\sqrt{v_{n,j}}}{2}\right)\P\left(N^-=\tilde\mu_{n,j}-\frac{\sqrt{v_{n,j}}}{2}\right)\\
&\geq & \frac{v_{n,j}}{(\ln
n)^{2\al}}\frac{\mu_{n,j}^{l_{n,j}}}{l_{n,j}!}e^{-\mu_{n,j}}\frac{\mu_{n,j}^{m_{n,j}}}{m_{n,j}!}e^{-\mu_{n,j}}\\
&\geq & \frac{v_{n,j}e^{-2}}{(\ln
n)^{2\al}}\frac{\tilde\mu_{n,j}^{l_{n,j}}}{l_{n,j}!}e^{-\tilde\mu_{n,j}}\frac{\tilde\mu_{n,j}^{m_{n,j}}}{m_{n,j}!}e^{-\tilde\mu_{n,j}}\\
&\geq &\frac{4\gamma'e^{-2}\tilde\mu_{n,j}}{(\ln
n)^{2\al-1}}\left(\frac{\tilde\mu_{n,j}}{l_{n,j}}\right)^{l_{n,j}}
e^{-(\tilde\mu_{n,j}-l_{n,j})}\left(\frac{\tilde\mu_{n,j}}{m_{n,j}}\right)^{m_{n,j}}
e^{-(\tilde\mu_{n,j}-m_{n,j})}\frac{(1+o_n(1))}{2\pi\sqrt{l_{n,j}m_{n,j}}}\\
&\geq & \frac{2\gamma'e^{-2}}{\pi(\ln
n)^{2\al-1}}e^{-\tilde\mu_{n,j} \left[h\left(\frac{\sqrt{v_{n,j}} }{2\tilde\mu_{n,j}}\right)+h\left(-\frac{\sqrt{v_{n,j}} }{2\tilde\mu_{n,j}}\right)\right]}(1+o_n(1))
\end{eqnarray*}
where $h(x)=(1+x)\ln(1+x)-x =x^2/2+O(x^3)$.
So, 
$$\E(\normp{\tilde{f}_{n,\gamma}^H-f}^2)\geq \frac{2\gamma'e^{-2}}{\pi(\ln
n)^{2\al-1}}e^{- \frac{v_{n,j}}{4\tilde\mu_{n,j}}+O_n\left(\frac{v_{n,j}^{\frac{3}{2}}}{\tilde\mu_{n,j}^2}\right)}(1+o_n(1)).$$
Since
$$v_{n,j}=4\gamma'\ln(n)\tilde\mu_{n,j}(1+o_n(1)),$$
we obtain
$$\E(\normp{\tilde{f}_{n,\gamma}^H-f}^2)\geq \frac{2\gamma'e^{-2}}{\pi(\ln
n)^{2\al-1}}e^{- \gamma'\ln(n)+o_n(\ln(n))}(1+o_n(1)).$$
Finally, for every $\delta>\gamma'$,
$$\E(\normp{\tilde{f}_{n,\gamma}^H-f}^2)\geq \frac{1}{n^{\delta}}(1+o_n(1)),$$
 and Theorem \ref{lower} is proved.
%%%%%%%%%%%%%%%%%%%%%%%%%%%%%%%%%%%%%%%%%%%%%%%%
\subsection{Proof of Theorem \ref{uppth}}
Without \pa{loss of} generality, the result is proved for $R=2$. Before proving Theorem \ref{uppth}, let us state the following result. 
\begin{Lemma}\label{gagrand}
Let  $\ga_{\min}\in(1,\gamma)$  be  fixed  and  let $\eta_{\la, \ga_{\min}}$  be  the  threshold
associated with $\ga_{\min}$:
 $$\eta_{\la, \ga_{\min}}=\sqrt{2\gamma_{\min}\ln n \tilde{V}_{\la,n}^{\min} }+\frac{\gamma_{\min}\ln n}{3n}\norm{\p_\la}_\infty,$$
where    $$\tilde{V}_{\la,n}^{\min}=\hat{V}_{\la,n}+\sqrt{2\gamma_{\min}     \ln    n    \hat{V}_{\la,n}
\frac{\norm{\p_\la}_\infty^2}{n^2}}+3                \gamma_{\min}                \ln
n\frac{\norm{\p_\la}_\infty^2}{n^2}$$
(see (\ref{defthresh})). Let $u=(u_n)_n$ be a sequence of positive
numbers and$$\La_u=\left\{\la\in\Gamma_n:\quad\P(\eta_{\la,\gamma} \leq|\be_\la|+\eta_{\la,\gamma_{\min}})\leq 
u_n\right\}.$$
Then
$$\E(\normp{\tilde{f}_{n,\gamma}^H-f}^2)\geq \left(\sum_{\la \in
\La_u}\be_\la^2\right) (1-(3n^{-\ga_{\min}}+u_n)).$$
\end{Lemma}
\begin{proof}
\begin{eqnarray*}
\E(\normp{\tilde{f}_{n,\gamma}^H-f}^2) &\geq & \sum_{\la \in
\La_u} \E\left((\hb_\la-\be_\la)^2\indic_{|\hb_\la|\geq \eta_{\la,\gamma}}+\be_\la^2 \indic_{|\hb_\la|<\eta_{\la,\gamma}}\right)\\
&\geq &\sum_{\la \in
\La_u}\be_\la^2 \P( |\hb_\la|<\eta_{\la,\gamma})\\
&\geq & \sum_{\la \in
\La_u}\be_\la^2 \P(|\hb_\la-\be_\la|+|\be_\la|<\eta_{\la,\gamma})\\
&\geq & \sum_{\la \in
\La_u}\be_\la^2 \P(|\hb_\la-\be_\la|<\eta_{\la,\gamma_{\min}} \mbox{ and }
\eta_{\la,\gamma_{\min}}+|\be_\la|<\eta_{\la,\gamma})\\
&\geq & \sum_{\la \in
\La_u}\be_\la^2
\left(1-\left(\P(|\hb_\la-\be_\la|\geq\eta_{\la,\gamma_{\min}})+\P(\eta_{\la,\gamma_{\min}}+|\be_\la|\geq\eta_{\la,\gamma})\right)\right)\\
&\geq &\left(\sum_{\la \in
\La_u}\be_\la^2\right) (1-(3n^{-\ga_{\min}}+u_n)),
\end{eqnarray*}
by applying \pa{the technical} Lemma \ref{toutesdev} of the Appendix section. 
\end{proof}
Using Lemma \ref{gagrand}, we give the proof of Theorem \ref{uppth}.
Let us consider 
$$f=\indic_{[0,1]} +\sum_{k\in\mathcal{N}_j}\sqrt{\frac{2(\sqrt{\ga}-\sqrt{\ga_{\min}})^2\ln n}{n}} \tp_{\jk},$$
with $$ \mathcal{N}_j=\{0,1,\dots,2^j-1\}$$
and   $$\frac{n}{(\ln  n)^{1+\al}}<2^j\leq   \frac{2n}{(\ln  n)^{1+\al}},\quad
\al>0.$$
Note that for any $k\in \mathcal{N}_j$, $$F_{j,k}=2^{-j}\geq \frac{(\ln n)(\ln\ln
n)}{n}$$ for $n$ large enough and $f$ belongs to ${\cal F}_n(2)$.  Furthermore, for any
$k\in \mathcal{N}_j$,
$$V_{(j,k),n}=V_{(-1,0),n}=\frac{1}{n}.$$ So, for $n$ large enough,
$$\sum_{\la \in \Ga_n}\min( \beta_\la^2, V_{\la,n})=V_{(-1,0),n}+\sum_{k\in\mathcal{N}_j}V_{(j,k),n}=\frac{1}{n}+\sum_{k\in\mathcal{N}_j}\frac{1}{n}.$$
Now,  to   apply  Lemma   \ref{gagrand},  let  us   set  for   any  $n$,
$u_n=n^{-\ga}$ and observe that for any $\e>0$, since $\ga_{\min}<\ga,$
$$
\P(\eta_{\la,\gamma_{\min}}+|\be_\la|\geq \eta_{\la,\gamma})\leq \P((1+\e)2\ga_{\min}\ln
n \tilde{V}_{\la,n}^{\min} + (1+\e^{-1}) \be_\la^2 > 2\ga\ln n
\tilde{V}_{\la,n}),
$$
with 
$$\be_\la^2=\frac{2(\sqrt{\ga}-\sqrt{\ga_{\min}})^2\ln n}{n}.$$
With $\e = \sqrt{\ga/\ga_{\min}}-1$ and $\theta=\sqrt{\ga_{\min}/\ga}$, 
$$\P((1+\e)2\ga_{\min}\ln
n \tilde{V}_{\la,n}^{\min} + (1+\e^{-1}) \be_\la^2 > 2\ga\ln n
\tilde{V}_{\la,n})= \P(\theta \tilde{V}_{\la,n}^{\min}+(1-\theta) V_{\la,n} > \tilde{V}_{\la,n}).$$
Since $ \tilde{V}_{\la,n}^{\min}<\tilde{V}_{\la,n}$, 
$$\P\left(\eta_{\la,\gamma_{\min}}+|\be_\la|\geq \eta_{\la,\gamma}\right)\leq \P(V_{\la,n} > \tilde{V}_{\la,n}) \leq u_n.$$
So, 
$$\{(j,k):\quad k\in\mathcal{N}_j\}\subset \La_u,$$
and 
\begin{eqnarray*}
\E(\normp{\tilde{f}_{n,\gamma}^H-f}^2)&\geq&\sum_{
k\in\mathcal{N}_j}\be_{j,k}^2(1-(3n^{-\ga_{\min}}+n^{-\ga}))\\
&\geq&(\sqrt{\ga}-\sqrt{\ga_{\min}})^2                                      2\ln
n\sum_{k\in\mathcal{N}_j}\frac{1}{n}(1-(3n^{-\ga_{\min}}+n^{-\ga}))\\
&\geq&(\sqrt{\ga}-\sqrt{\ga_{\min}})^2 2\ln n\left(\sum_{\la \in \Ga_n}\min( \beta_\la^2, V_{\la,n})-\frac{1}{n}\right)(1-(3n^{-\ga_{\min}}+n^{-\ga})).
\end{eqnarray*}
Finally, since $\mbox{card}(\mathcal{N}_j)\to +\infty$ when $n\to +\infty$,
$$\frac{\E(\normp{\tilde{f}_{n,\gamma}-f}^2)}{\sum_{\la \in \Ga_n}
\min( \beta_\la^2, V_{\la,n})+\frac{1}{n}}\geq (\sqrt{\ga}-\sqrt{\ga_{\min}})^2 2\ln n(1+o_n(1)).$$
%%%%%%%%%%%%%%%%%%%%%%%%%%%%%%%%%%%%%%%%%%%
%%%%%%%%%%%%%%%%%%%%%%%%%%%%%%%%%%%%%%%%%%%
\section{Appendix: Technical tools}\label{maintools}
%%%%%%%%%%%%%%%%%%%%%%%%%%%%%%%%%%%%%%%%%%%

%%%%%%%%%%%%%%%%%%%%%%%%%%%%
\subsection{Some probabilistic properties of the Poisson process}\label{rappelPoisson}
Let us first recall some basic facts about Poisson processes. 
\begin{Def}\label{poisson}
Let $(X,\Xx)$ be a measurable space.  Let $N$ be a random countable subset of
$X$. $N$ is said to be a Poisson process on $(X,\Xx)$ if
\begin{enumerate}
\item for any $A \in\Xx$, the number of points of $N$ lying in $A$ is a random
variable, denoted $N_A$, which obeys a Poisson distribution with parameter $\mu(A)$,
where $\mu$ is a measure on $X$.
\item for any finite family of disjoint sets $A_1,...,A_n$ of $\Xx$, $\N{A_1},...,\N{A_n}$ are independent random variables. 
\end{enumerate}
\end{Def}
\pa{We focus here on the case } $X=\R$. \pa{Let us } mention that a Poisson process $N$ is infinitely
 divisible, which means that it can be written as follows: for any positive integer $k$:
 \begin{equation}  \label{infdiv}
   dN=\sum\limits_{i=1}^k~ dN_i
 \end{equation}
where the $N_i$'s are mutually %({\bf on peut enlever mutually ?}) 
independent Poisson processes on
 $\R$ with mean measure $\mu/k$. The following
 proposition (sometimes attributed to  Campbell (see \cite{kin})) is fundamental.
 \begin{Prop}
\label{esperances}
 For any measurable function $g$ and  any $z\in\R$, such that $\int
e^{zg(x)}d\mu_x<\infty$ one has,
$$\E\left[\exp\left(z\int_\R g(x)
dN_x\right)\right]=\exp\left(\int_\R \left(e^{z g(x)}-1\right)
d\mu_x\right).$$ So,
$$
\E\left(\int_\R g(x) dN_x\right)=\int_\R g(x) d\mu_x,\quad
\var\left(\int_\R g(x) dN_x\right)=\int_\R g^2(x) d\mu_x.
$$
If $g$ is bounded, this implies the following
exponential inequality. For any $u>0$,
\begin{equation}
 \label{invformben}
\P\left(\int_\R g(x) (dN_x-d\mu_x) \geq
   \sqrt{2u\int_\R g^2(x) d\mu_x}+\frac{1}{3}\norm{g}_\infty
    u\right)\leq\exp(-u).
\end{equation}
\end{Prop}
%%%%%%%%%%%%%%%%%%%%%%%%%%%%
\subsection{Biorthogonal wavelet bases}\label{biorthogonal}
We set 
$$\phi=\indic_{[0,1]}.$$ For any $r> 0$, there exist three functions
$\psi$, $\tilde \phi$ and $\tilde\psi$ with the following properties:
\begin{enumerate}
\item $\tilde\phi$ and $\tilde\psi$ are compactly supported,
\item  $\tilde\phi$  and $\tilde\psi$  belong  to  $C^{r+1}$, where  $C^{r+1}$
denotes the H\"older space of order $r+1$,
\item $ \psi$ is compactly supported and is a piecewise constant function,
\item $\psi$ is orthogonal to polynomials of degree no larger than $r$,
\item $\{(\phi_{k},\psi_{j,k})_{j\geq
0,k\in\Z},(\tilde\phi_{k},\tilde\psi_{j,k})_{j\geq     0,k\in\Z}\}$    is    a
biorthogonal family: for any $j,j'\geq 0,$ for any $k,k',$
$$\int_\R\psi_{j,k}(x)\tilde\phi_{k'}(x)dx=\int_\R\phi_{k}(x)\tilde\psi_{j',k'}(x)dx=0,$$
$$\int_\R\phi_{k}(x)\tilde\phi_{k'}(x)dx=1_{k=k'},\quad
\int_\R\psi_{j,k}(x)\tilde\psi_{j',k'}(x)dx=1_{j=j',k=k'},$$
where for any $x\in\R$ and for any $(j,k)\in~\Z^2$,
$$\phi_{k}(x)=\phi(x-k), \quad \psi_{j,k}(x)=2^{\frac{j}{2}}\psi(2^jx-k)$$
and
$$\tilde\phi_{k}(x)=\tilde\phi(x-k), \quad \tilde\psi_{j,k}(x)=2^{\frac{j}{2}}\tilde\psi(2^jx-k).$$
\end{enumerate}
This implies the wavelet decomposition (\ref{decom2}) of $f$. Such biorthogonal wavelet bases have been built by
Cohen Daubechies and Feauveau \cite{cdf} as a special
case  of spline  systems (see  also the  elegant  equivalent construction  of
Donoho \cite{don}  from boxcar functions). The  Haar basis can be  viewed as a
particular biorthogonal wavelet basis, by setting  $\tilde\phi=\phi$ and
$\tilde\psi=\psi=\indic_{[0,\frac{1}{2}]}-\indic_{]\frac{1}{2},1]}$, with $r=0$  even if Property 2
is not satisfied with such a choice. The Haar basis is 
%Of course, recall  that all these properties except the  second and the forth
%ones   are   true  for   the   Haar   basis,   where  $\tilde\phi=\phi$   and
%$\tilde\psi=\psi=\indic_{[0,1/2]}-\indic_{]1/2,1]}$,  which allows  to obtain
%in addition 
an orthonormal basis \pa{but this} is not true for general biorthogonal wavelet
bases. However, we have the frame property: if we denote
$$\varPhi=\{\phi,\psi,\tilde\phi,\tilde\psi\}$$
there exist two constants $c_1(\varPhi)$ and $c_2(\varPhi)$ only depending on $\varPhi$ such that
$$c_1(\varPhi)\left(\sum_{k\in\Z}\alpha_k^2+\sum_{j\geq
0}\sum_{k\in\Z}\beta_{j,k}^2\right)\leq                         \|f\|_{2}^2\leq
c_2(\varPhi)\left(\sum_{k\in\Z}\alpha_k^2+\sum_{j\geq
0}\sum_{k\in\Z}\beta_{j,k}^2\right).$$
For    instance,    when    the    Haar    basis    is    considered,
$c_1(\varPhi)=c_2(\varPhi)=1$. In particular, we have
\begin{equation}\label{equnormes}
c_1(\varPhi)\norm{\tilde\be-\be}_{\ell_2}^2\leq \|\tilde f_{n,\ga}-f\|_{2}^2\leq
c_2(\varPhi)\norm{\tilde\be-\be}_{\ell_2}^2. 
\end{equation}
 An important feature of such bases  is the following: there exists a constant $\mu_{\psi}>0$ such that
\begin{equation}\label{minophi}
\inf_{x\in[0,1]}|\phi(x)|\geq 1,\quad\inf_{x\in \supp(\psi)}|\psi(x)|\geq\mu_{\psi},
\end{equation}
where $\supp(\psi)=\{x\in\R:\quad \psi(x)\not=0\}.$
%%%%%%%%%%%%%%%%%%%%%%%%%%%%%%%%%%%%%%%%%%%
\subsection{Proof of Lemma \ref{lemmesharp}}\label{preuvelemmesharp}
The proof of Lemma \ref{lemmesharp} is based on the following result proved in \cite{Poisson_minimax}.
\begin{Th}
\label{inegmodelsel}
To estimate a countable family $\be=(\be_\la)_{\la\in \La}$, such that $\|\be\|_{\ell_2}<\infty$, we assume that
 a family of coefficient estimators $(\hb_\la)_{\la \in \Ga}$, where $\Ga$ is
a known deterministic subset of $\La$, and a family of possibly
random thresholds $(\eta_\la)_{\la\in \Ga}$ are available. \pa{We} consider
the thresholding rule $\tb=(\hb_\la\indic_{|\hb_\la|\geq
  \eta_\la}\indic_{\la\in\Ga})_{\la \in \La}$.
Let $\e>0$ be fixed.
Assume  that there exist  a deterministic  family $(F_\la)_{\la\in\Ga}$ and
three constants $\kappa\in [0,1[$,  $\omega\in [0,1]$ and $\mu>0$ (that
may depend on $\e$ but not on $\la$) with the following properties. 
\begin{itemize}
\item[(A1)] For   all  $\la$  in   $\Ga$,  $$\P(|\hb_\la-\be_\la|>\kappa\eta_\la)\leq\omega.$$
\item[(A2)] There exist $1<p,q<\infty$ with $\frac{1}{p}+\frac{1}{q}=1$ and a constant $R>0$ such that for all $\la$ in $\Ga$,
  $$\left(\E(|\hb_\la-\be_\la|^{2p})\right)^{\frac{1}{p}}\leq R \max(F_\la,
  F_\la^{\frac{1}{p}}\e^{\frac{1}{q}}).$$
\item[(A3)] There exists  a constant $\theta$ such that for all  $\la$ in $\Ga$ such
that $F_\la<\theta \e$
$$\P(|\hb_\la-\be_\la|>\kappa\eta_\la, |\hb_\la|>\eta_\la)\leq F_\la\mu.$$ 
\end{itemize}
Then the estimator $\tb$ satisfies
$$\frac{1-\kappa^2}{1+\kappa^2}\E\|\tb-\be\|_{\ell_2}^2\leq \E\inf_{m\subset
  \Ga}\left\{\frac{1+\kappa^2}{1-\kappa^2}\sum_{\la\not \in
  m}\be_\la^2+\frac{1-\kappa^2}{\kappa^2}\sum_{\la\in
  m}(\hb_\la-\be_\la)^2+\sum_{\la \in m}\eta_\la^2\right\}+LD\sum_{\la\in\Ga}F_\la
$$
with $$LD=\frac{R}{\kappa^2}\left(\left(1+\theta^{-1/q}\right)\omega^{1/q}+(1+\theta^{1/q})\e^{1/q}\mu^{1/q}\right).$$
\end{Th}
To prove Lemma \ref{lemmesharp},  we apply Theorem \ref{inegmodelsel} with $\hb_\la$  defined in  (\ref{defest1}),   $\eta_\la=\eta_{\la,\ga}$  defined   in
(\ref{defthresh})  and 
$\Ga=\Ga_n$ defined in (\ref{Gamman}). 
We set
$$F_\la=\int_{\supp(\p_\la)}f(x)dx,$$
so we have:
\begin{equation}\label{Fla}
\sum_{\la\in\Ga_n}F_\la=\sum_{-1\leq j\leq j_0}\sum_k\int_{x\in
\supp(\p_{j,k})} f(x)dx\leq\int f(x)dx\sum_{-1\leq j\leq
j_0}\sum_k\indic_{x\in \supp(\p_{j,k})} \leq (j_0+2)m_\p\norm{f}_1,
\end{equation}
where $m_\p$  is a finite constant depending  only on the  compactly  supported functions $\phi$ and
$\psi$.   Finally, $\sum_{\la\in\Ga_n}F_\la$  is
bounded by $\log(n)$ up to a  constant that only depends on $\norm{f}_1$ and the functions $\phi$ and $\psi$. Now, we give a fundamental lemma to
derive Assumption (A1) of  Theorem \ref{inegmodelsel}.
\begin{Lemma} 
\label{toutesdev}
For any $u>0$,
\begin{equation}
\label{surhb1}
\P\left(|\hb_\la-\be_\la|\geq
\sqrt{2uV_{\la,n}}+\frac{\norm{\p_\la}_\infty u}{3n}\right)\leq 2 e^{-u}.
\end{equation}
Moreover, for any $u>0$,
\begin{equation}
\label{surv1}
\P\left(V_{\la,n}\geq \breve{V}_{\la,n}(u)\right)\leq e^{-u},\nonumber
\end{equation}
where
$$\breve{V}_{\la,n}(u) = \hat{V}_{\la,n} +\sqrt{2\hat{V}_{\la,n} \frac{\norm{\p_\la}_\infty^2}{n^2}u}+3\frac{\norm{\p_\la}_\infty^2}{n^2}u.$$
\end{Lemma}
\begin{proof}
Equation  (\ref{surhb1})  comes easily  from  (\ref{invformben}) applied  with
$g=\p_\la/n$. The same inequality applied with $g=-\p_\la^2/n^2$ gives:
$$\P\left(V_{\la,n}\geq \hat{V}_{\la,n} +\sqrt{2u \int_\X \frac{\p_\la^4(x)}{n^4}
n f(x) dx}+\frac{\norm{\p_\la}_\infty^2}{3n^2}u \right)\leq e^{-u}.$$
We observe that $$\int_\X \frac{\p_\la^4(x)}{n^4}
n f(x) dx\leq \frac{\norm{\p_\la}_\infty^2}{n^2} V_{\la,n}.$$
So, if we set $a=u\frac{\norm{\p_\la}_\infty^2}{n^2}$, then
$$\P(V_{\la,n}-\sqrt{2V_{\la,n} a}-a/3\geq \hat{V}_{\la,n})\leq e^{-u}.$$
%Let $\mathcal{P}(x)=x^2-\sqrt{2a}x-a/3$. The discriminant of this polynomial is $10a/3$ which is strictly larger than $2a$. Since $V_{\la,n}$ and $\hat{V}_{\la,n}$ are positive, 
We obtain
$$\P(\sqrt{V_{\la,n}}\geq \mathcal{P}^{-1}(\hat{V}_{\la,n}))\leq e^{-u}$$
where $\mathcal{P}^{-1}(\hat{V}_{\la,n})$ is the positive solution of 
$$(\mathcal{P}^{-1}(\hat{V}_{\la,n}))^2-\sqrt{2a}\mathcal{P}^{-1}(\hat{V}_{\la,n})-(a/3+\hat{V}_{\la,n})=0.$$
To conclude, it remains to observe that $$\breve{V}_{\la,n}(u) \geq(\mathcal{P}^{-1}(\hat{V}_{\la,n}))^2 =\left(\sqrt{\hat{V}_{\la,n}+5a/6}+\sqrt{a/2}\right)^2.$$
\end{proof}
Let $\kappa<1$. Combining these inequalities with $\tilde{V}_{\la,n}=\breve{V}_{\la,n}(\ga \ln n )$ yields 
\begin{eqnarray*}
 \P(|\hb_\la-\be_\la|>\kappa\eta_{\la,\gamma}) 
& \leq& \P\left(|\hb_\la-\be_\la|\geq \sqrt{2\kappa^2\ga\ln n 
    \tilde{V}_{\la,n}}+\frac{\kappa\gamma\ln n
    \norm{\p_\la}_\infty }{3 n}\right)\\
&\leq& \P\left(|\hb_\la-\be_\la|\geq\sqrt{2\kappa^2\ga\ln n 
    \tilde{V}_{\la,n}}+\frac{\kappa\gamma\ln n
    \norm{\p_\la}_\infty }{3 n} , V_{\la,n} \geq \tilde{V}_{\la,n}\right)\\&&+\P\left(|\hb_\la-\be_\la|\geq\sqrt{2\kappa^2\ga\ln n 
    \tilde{V}_{\la,n}}+\frac{\kappa\gamma\ln n
    \norm{\p_\la}_\infty }{3 n}, V_{\la,n} < \tilde{V}_{\la,n}\right)\\
&\leq &\P(V_{\la,n} \geq \tilde{V}_{\la,n})+ \P\left(|\hb_\la-\be_\la|\geq \sqrt{2\kappa^2\ga\ln n V_{\la,n}}+\frac{\kappa\gamma\ln n \norm{\p_\la}_\infty }{3n}\right)\\
&\leq& n^{-\ga}+2n^{-\kappa^2\ga}\\&\leq &3 n^{-\kappa^2\ga}.
\end{eqnarray*}
So, for any value of $\kappa\in [0,1[$, Assumption (A1) is true with $\eta_{\lambda}=\eta_{\lambda,\gamma}$ and $\Gamma=\Gamma_n$ if we take $\omega=3 n^{-\kappa^2\ga}$. To \pa{satisfy} the Rosenthal type inequality (A2) of Theorem \ref{inegmodelsel}, we
prove the following lemma.
\begin{Lemma}
\label{moment2}
For any $p\pa{>1}$, there exists an absolute constant $C$ such that
$$\E(|\hb_\la-\be_\la|^{2p}) \leq C^p p^{2p} \left(V_{\la,n}^p+\left[\frac{\norm{\p_\la}_\infty}{n}\right]^{2p-2}V_{\la,n}\right).$$
\end{Lemma}
\begin{proof}
We apply (\ref{infdiv}). Hence, 
$$\hb_\la-\be_\la=\sum_{i=1}^k \int \frac{\p_\la(x)}{n} \left(dN_x^i-nk^{-1}f(x)dx\right)=\sum_{i=1}^k Y_i$$
where for any $i$, $$Y_i= \int \frac{\p_\la(x)}{n} \left(dN_x^i-nk^{-1}f(x)dx\right).$$ So
the $Y_i$'s  are i.i.d. centered  variables, each of  them \pa{having} a moment  of order
$2p$.  For any $i$, we apply the Rosenthal inequality (see Theorem 2.5 of \cite{johnson}) to
the positive and negative parts of $Y_i$. This easily implies that
$$\E\left(\left|\sum_{i=1}^k   Y_i\right|^{2p}\right)    \leq   \left(\frac{16p}{\ln(2p)}\right)^{2p}   \max\left(\left(\E\sum_{i=1}^kY_i^2\right)^p,
\left(\E\sum_{i=1}^k |Y_i|^{2p}\right)\right).$$ 
It remains to bound the upper limit of $\E(\sum_{i=1}^k |Y_i|^{\ell})$ for all
$\ell\in\{2p,2\}\geq     2$    when     $k\to\infty$.     Let    us     introduce
$$\Omega_k=\{ \mbox{card}(N^i_{\pa{\R}})\leq 1 \ \mbox{for any } i\in\{1,\dots,k\}\}.$$ Then,  it is easy
to see that $\P(\Omega_k^c)\leq k^{-1}(n\norm{f}_1)^2 $ (see e.g., (\ref{n2}) below). 
\\\\
On $\Omega_k$, $|Y_i|^\ell= O_k(k^{-\ell})$ if \pa{$\mbox{card}(N^i_{\pa{\R}})=0$ }%$\int \frac{\p_\la(x)}{n} dN_x^i =0$ 
and $|Y_i|^\ell= \left[\frac{|\p_\la(T)|}{n}\right]^{\ell} + O_k\left(k^{-1}\left[\frac{|\p_\la(T)|}{n}\right]^{\ell-1}\right)$ if $\int \frac{\p_\la(x)}{n} dN_x^i =\frac{\p_\la(T)}{n} $ where $T$ is the point of the process $N^i$.
Consequently,
\begin{multline}
\label{split}
\E\sum_{i=1}^k |Y_i|^{\ell} \leq \E\left(\indic_{\Omega_k} \left(\sum_{T\in N}\left[\left[\frac{|\p_\la(T)|}{n}\right]^{\ell} + O_k\left(k^{-1}\left[\frac{|\p_\la(T)|}{n}\right]^{\ell-1}\right)\right] + k O_k(k^{-\ell})\right)\right)\\ + \sqrt{\P(\Omega_k^c)} \sqrt{\E\left[\left(\sum_{i=1}^k |Y_i|^\ell\right)^2\right]}.
\end{multline}
But we have\begin{eqnarray*}\sum_{i=1}^k |Y_i|^\ell &\leq& 2^{\ell-1} \left(\sum_{i=1}^k \left[\left[\frac{\norm{\p_\la}_\infty}{n}\right]^\ell (N^i_\X)^\ell + \left(k^{-1}\int |\p_\la(x)| f(x) dx \right)^\ell\right]\right)\\
&\leq& 2^{\ell-1} \left(\left[\frac{\norm{\p_\la}_\infty}{n}\right]^\ell (N_\X)^\ell + k\left(k^{-1}\int |\p_\la(x)| f(x) dx \right)^\ell\right).\end{eqnarray*}
So, when $k\to +\infty$, the last term in (\ref{split}) converges to 0 since a Poisson variable has moments of every order and
$$
\mbox{limsup}_{k\to\infty} \E\sum_{i=1}^k |Y_i|^{\ell} \leq  \E\left( \int \left[\frac{|\p_\la(x)|}{n}\right]^{\ell} dN_x\right)\leq \left[\frac{\norm{\p_\la}_\infty}{n}\right]^{\ell-2} V_{\la,n},$$ 
which concludes the proof.
\end{proof}
Now,
\begin{equation}\label{relFV}
V_{\la,n}=\frac{1}{n}\int\p_\la^2(x)f(x)dx\leq \frac{\norm{\p_\la}_\infty^2 F_\la}{n} 
\end{equation}
and    Assumption    (A2)    is    satisfied   with    $\e=\frac{1}{n}$    and
$$R=\frac{2Cp^22^{j_0}\max(\norm{\phi}_\infty^2;\norm{\psi}_\infty^2)}{n}$$
since                                                   $\norm{\p_\la}_\infty^2\leq
2^{j_0}\max(\norm{\phi}_\infty^2;\norm{\psi}_\infty^2)$ and
\begin{eqnarray*}
\left(\E(|\hb_\la-\be_\la|^{2p})\right)^{\frac{1}{p}}&\leq&Cp^2\left(\frac{\norm{\p_\la}_\infty^2
F_\la}{n}                                               +\norm{\p_\la}_\infty^2
F_\la^{\frac{1}{p}}n^{\frac{1}{p}-2}\right)
\leq\frac{Cp^2\norm{\p_\la}_\infty^2}{n}\left(F_\la+F_\la^{\frac{1}{p}}n^{-\frac{1}{q}}\right).
\end{eqnarray*}
Finally, Assumption (A3) comes from the following lemma.
\begin{Lemma}
\label{nombredepoints} We set $$N_\la=\int_{\supp(\p_\la)}dN \quad\mbox{ and }\quad C'=(\sqrt{6}+1/3)\gamma\geq\sqrt{6}+1/3.$$ There exists an absolute constant $0<\theta'<1$ such
that if 
$$
nF_\la \leq \theta' C' \ln n
$$
 and 
\begin{equation}\label{nF}
(1-\theta')(\sqrt{6}+1/3) \ln n \geq 2\end{equation} then,
$$\P(N_\la-nF_\la \geq (1-\theta') C' \ln n) \leq F_\la n^{-\gamma}.$$
\end{Lemma}
\begin{Remark}\label{nvalide}
We can take $\theta'=0.01$ and in this case, (\ref{nF})  is satisfied as soon as $n\geq 3$. 
\end{Remark}
\begin{proof}
One takes $\theta'\in [0,1]$ (for instance $\theta'=0.01$) such that $$\frac{3(1-\theta')^2}{2(2\theta'+1)}(\sqrt{6}+1/3)\geq 4.$$
We use Equation (5.2) of \cite{ptrfpois} to obtain
$$\P(N_\la-nF_\la \geq (1-\theta')C' \ln n)
\leq\exp\left(-\frac{((1-\theta')C'\ln n)^2}{2(nF_\la+(1-\theta')C'\ln n/3)}\right)\leq
n^{-\frac{3(1-\theta')^2}{2(2\theta'+1)}C'}.$$
If $nF_\la \geq n^{-\gamma-1}$, since  $\frac{3(1-\theta')^2}{2(2\theta'+1)} C' \geq
2 \gamma+2,$ the result is true. 
If  $nF_\la \leq n^{-\gamma-1}$, 
\begin{equation}\label{n2}
\P(N_\la-nF_\la \geq (1-\theta')C' \ln n)\leq \P(N_\la >(1-\theta')C'
\ln n)\leq \P(N_\la \geq 2)
\leq \sum_{k\geq 2} \frac{(nF_\la)^k}{k!} e^{-n F_\la}\leq (nF_\la)^2
\end{equation}
and the result is true.
\end{proof} 
Now, observe that if $|\hb_\la|>
\eta_{\la,\gamma}$ then 
$$N_\la\geq C' \ln n.$$ Indeed, $|\hb_\la|>
\eta_{\la,\gamma}$ implies
$$\frac{C'\ln n
}n \norm{\p_\la}_\infty\leq |\hb_\la| \leq \frac{\norm{\p_\la}_\infty
N_\la}{n}.$$
So if $n$ satisfies  $(1-\theta')(\sqrt{6}+1/3) \ln n \geq 2$, we set $\theta=\theta' C'\ln(n)$ and
$\mu=n^{-\gamma}$. In this case, Assumption (A3) is fulfilled since if $nF_\la \leq \theta' C' \ln n $ 
$$\P(|\hb_\la-\be_\la|>\kappa\eta_\la, |\hb_\la|>\eta_\la)\leq \P(N_\la-nF_\la
\geq (1-\theta') C' \ln n) \leq F_\la n^{-\gamma}.$$
Finally, if $n$ satisfies  $(1-\theta')(\sqrt{6}+1/3) \ln n \geq 2$, Theorem \ref{inegmodelsel} gives:
$$
\frac{1-\kappa^2}{1+\kappa^2}\E\norm{\tb-\be}^2_{\ell_2}\leq \inf_{m\subset
  \Ga_{\pa{n}}}\left\{\frac{1+\kappa^2}{1-\kappa^2}\sum_{\la\not \in
  m}\be_\la^2+\frac{1-\kappa^2}{\kappa^2}\sum_{\la\in
  m}\E(\hb_\la-\be_\la)^2+\sum_{\la \in m}\E(\eta_{\la,\gamma}^2)\right\}+LD\sum_{\la\in\Ga}F_\la.
$$
In addition,  there exists a constant  $K_1$ depending on  $p$, $\gamma$, $\norm{f}_1$ and on $\Phi$ such that 
\begin{equation}\label{init2}
LD\sum_{\la\in\Ga}F_\la\leq K_1\log(n)n^{-\frac{\kappa^2\ga}{q}}.
\end{equation} 
Since $\ga>1$, \pa{for all} $\kappa<1$, \pa{there exists} $q>1$ such that $1<\frac{\kappa^2\ga}{q}$ and as required by Theorem \ref{inegoraclelavraie}, the last term satisfies
$$LD\sum_{\la\in\Ga}F_\la\leq \frac{\pa{K(\gamma,\kappa,\norm{f}_1)}}{n},$$
\pa{where $K(\gamma,\kappa,\norm{f}_1)$ denotes a positive constant. This concludes the proofs.}
%%%%%%%%%%%%%%%%%%%%%%%%%%%%%%%%%%%%%%%%%%%%%%%%%%%%%%%%%%%%%%%%%
%%%%%%%%%%%%%%%%%%%%%%%%%%%%%%%%%%%%%%%%%%%%%%%%%%%%%%%%%%%%%%%%%
\section{Definition of the signals used in Section \ref{simulations}}\label{defsignal}
\noindent The following table gives the definition of the signals used in Section \ref{simulations}.

\vspace{1cm}
{\scriptsize
\hspace{-0.8cm}
\begin{tabular}{|c|c|c|}
\hline 
&&\\
Haar1 & Haar2 &  Blocks \\ 
&&\\
 $\displaystyle{\bf 1}_{[0,1]}$ & $\displaystyle  1.5~{\bf 1}_{[0,0.125]}+0.5~{\bf 1}_{[0.125,0.25]}+{\bf 1}_{[0.25,1]}$ & $\displaystyle\left(2+\sum_j \frac{h_j}2\left(1+\mbox{sgn}(x-p_j)\right)\right)\frac{{\bf 1}_{[0,1]}}{3.551}$ \\ 
&&\\
\hline
&&\\
Comb  & Gauss1 & Gauss2  \\ 
&&\\
$\displaystyle 32\sum_{k=1}^{+\infty}\frac1{k2^k}{\bf
  1}_{[k^2/32,(k^2+k)/32]}$& $\displaystyle\frac{1}{0.25\sqrt{2\pi}}
\exp\left(\frac{-(x-0.5)^2}{2\times 0.25^2}\right)$ &$\displaystyle
\frac{1}{\sqrt{2\pi}}
\exp\left({\frac{-(x-0.5)^2}{2\times 0.25^2}}\right)+\frac{3}{\sqrt{2\pi}}
\exp\left({\frac{-(x-5)^2}{2\times 0.25^2}}\right)$\\
&&\\
\hline 
&&\\
Beta0.5 & Beta4 & Bumps\\
&&\\
$\displaystyle0.5x^{-0.5}{\bf 1}_{]0,1]}$ & $\displaystyle 3x^{-4}{\bf 1}_{[1,+\infty[}$ & $\displaystyle\left(\sum_jg_j\left(1+\frac{|x-p_j|}{w_j}\right)^{-4}\right)\frac{{\bf 1}_{[0,1]}}{0.284}$\\
&&\\
\hline
\end{tabular}}
\vspace{0.3cm}

\noindent where

$ 
\begin{tabular}{ccccccccccccccc}
  p & =  & [ & 0.1 & 0.13 & 0.15 &0.23 &0.25& 0.4  & 0.44 & 0.65 & 0.76 &0.78  & 0.81 &  ]  \\ 
  h & =  & [  & 4 & -5 &3  & -4 & 5 &-4.2  &2.1  & 4.3 & -3.1 & 2.1 & -4.2 & ]  \\ 
 g  & = & [ & 4  & 5 & 3 & 4 & 5 & 4.2 & 2.1 & 4.3  &3.1  & 5.1  & 4.2  & ] \\ 
  w&=  &[  &0.005  & 0.005 &0.006  &0.01  &0.01  &0.03  &0.01  &0.01  &0.005  &0.008  &0.005  & ] 
\end{tabular} 
$
\\\\
\noindent{\bf Acknowledgment.} The authors  acknowledge the support of the  French Agence Nationale
de la Recherche  (ANR), under grant ATLAS (JCJC06\_137446) ''From Applications
to Theory in Learning and Adaptive Statistics''. We would like to warmly thank \pa{Rebecca Willett for her remarkable program, called FREE-DEGREE.}
%%%%%%%%%%%%%%%%%%%%%%%%%%%%%%%%%%%%%%%%%%%
%%%%%%%%%%%%%%%%%%%%%%%%%%%%%%%%%%%%%%%%%%%
\bibliographystyle{plain}

\begin{thebibliography}{10}
%\bibitem{abs}  Antoniadis,  A.,  Besbeas,   P.,  Sapatinas,  E.  {\it  Wavelet
%shrinkage  for natural  exponential families  with cubic  variance functions},
%Sankhya {\bf 63}, 309--327, (2001). 
%\bibitem{as}  Antoniadis,  A.,  Sapatinas,  T.  {\it  Wavelet
%shrinkage for natural exponential families with quadratic variance functions},
%Biometrika {\bf 88}(3), 805--820, (2001).
\bibitem{all} Allen, D.M. (1974).  The relationship between variable selection and data augmentation and a method for prediction. {\it Technometrics} {\bf 16} 125--127.
%\bibitem{antoniadisnew} Antoniadis, A., Sardy, S. and Tseng, P. (2004). Automatic smoothing with wavelets for a wide class of distributions.   {\it Journal of Computational and Graphical Statistics} {\bf 13}(2) 399--421.
\bibitem{am} Arlot, S. and  Massart, P. (2009). Data-driven calibration of penalties for least-squares regression. {\it Journal of Machine Learning Research} {\bf 10} 245--279.
\bibitem{ans} Anscombe, F.J. (1948). The transformation of Poisson, binomial and negative binomial data. {\it Biometrika} {\bf 35} 246--254.
\bibitem{aut} Autin, F. (2006). Maxiset for density estimation on $\mathbb R$.
{\it Math. Methods Statist} {\bf 15}(2) 123--145.
%\bibitem{apr} Autin F., Picard D., Rivoirard V. {\it Large variance Gaussian priors in Bayesian nonparametric estimation: a maxiset approach}, Mathematical Methods of Statistics, {\bf 15}(4), 349-373, (2006).
\bibitem{bb} Baraud, Y. and Birg\'e, L. (2006). Estimating the intensity of a
random measure by histogram type estimators. Technical report. To appear in Probab. Theory Related Fields.
\bibitem{bgh} Baraud, Y., Giraud, C. and Huet, S. (2008). Gaussian model selection with unknown variance. Technical report. To appear in Annals of Statistics.
%\bibitem{br} Bertin, K., Rivoirard, V. {\it Maxiset in sup-norm for kernel estimators}, to appear in Test,  (2008) 
\bibitem{bfs} Besbeas, P., De Feis, I. and Sapatinas, T. (2002). A Comparative Simulation Study of Wavelet Shrinkage Estimators For Poisson Counts. Technical report.
%\bibitem{birFano} Birg\'e, L. {\it A new look at an old result: Fano's Lemma}, 2001, manuscript. 
\bibitem{bir} Birg\'e, L. (2006). Model selection for Poisson processes. Technical report.
\bibitem{bm} Birg\'e, L. and Massart, P. (2007). Minimal penalties for
Gaussian model selection. {\it  Probab. Theory Related Fields} {\bf 138}(1-2) 33--73.
%\bibitem{bh} Bretagnolle,  J., Huber, C. {\it Estimation  des densit\'es: risque minimax}, Z. Wahrsch. Verw. Gebiete {\bf 47}(2), 119--137, (1979). 
%\bibitem{btw}  Bunea F.,  Tsybakov, A.B.,  Wegkamp, M.H.  {\it  Sparse density
%estimation with $l_1$ penalties}, 2007, manuscript. 
\bibitem{ck} Cavalier, L. and Koo, J.Y. (2002). Poisson intensity estimation for tomographic data using a wavelet shrinkage approach. 
{\it IEEE Trans. Inform. Theory} {\bf 48}(10) 2794--2802.
\bibitem{cdf} Cohen,  A., Daubechies,  I. and Feauveau, J.C. (1992). Biorthogonal
bases of  compactly supported wavelets. {\it  Comm. Pure Appl. Math.}   {\bf 45}(5)
485--560.
%\bibitem{coronel} Coronel-Brizio, H.F., Hernandez-Montoya, A.R. {\it On fitting the Pareto-L\'evy distribution to stock market index data: Selecting a suitable cutoff value} Physica A: Statistical Mechanics and its Applications, {\bf 354}, 437-449, (2005).
%\bibitem{del} Delyon, B. {\it Ondelettes orthogonales et biorthogonales}, 1996, manuscript. 
%\bibitem{dj}  Delyon, B.,  Juditsky, A.  {\it  On the  computation of  wavelet coefficients}, J. Approx. Theory {\bf 88}(1), 47--79, (1997). 
%\bibitem{dl} DeVore, R.A., Lorentz, G.G. {\it Constructive approximation}, Springer-Verlag, Berlin, 1993. 
%\bibitem{don2}  Donoho, D.L. {\it  Nonlinear wavelet  methods for  recovery of signals, densities, and spectra from indirect and noisy data},  Different perspectives on wavelets (San Antonio, TX, 1993),  173--205,Proc. Sympos. Appl. Math., {\bf 47}, Amer. Math. Soc., Providence, RI, (1993).
\bibitem{don}  Donoho, D.L. (1994).  Smooth wavelet  decompositions with  blocky
coefficient   kernels. {\it   Recent  advances   in   wavelet  analysis,   Wavelet
Anal. Appl.} {\bf 3} Academic Press, Boston, MA 259--308. 
\bibitem{dojo}  Donoho, D.L. and Johnstone,  I.M. (1994). Ideal spatial  adaptation by
wavelet shrinkage.  {\it  Biometrika} {\bf 81}(3) 425--455.
\bibitem{djkp}  Donoho,   D.L.,  Johnstone,  I.M.,   Kerkyacharian  G. and  Picard
D.  (1996).  Density estimation by wavelet thresholding. {\it Annals of Statistics}
{\bf 24}(2) 508--539.
%\bibitem{fh} Figueroa-L\'opez, J.E., Houdr\'e, C. {\it Risk bounds for the non-parametric estimation of L\'evy processes}, IMS Lecture Notes-Monograph series High Dimensional Probability, {\bf 51}, 96-116, (2006).
\bibitem{gei} Geisser, S.  (1975). The predictive sample reuse method with applications.  {\it J. Amer. Statist. Assoc.} {\bf 70} 320--328.
%\bibitem{gol} Golubev, G.K. {\it  Nonparametric estimation of smooth densities of a  distribution in  $L\sb 2$}, Problems  Inform. Transmission  {\bf 28}(1), 44--54, (1992). 
%\bibitem{gs} Gusto, G., Schbath, S.  {\it FADO: a statistical method to detect
%favored or avoided distances between motif occurrences using the hawkes model}, Statistical Applications in Genetics and Molecular Biology, {\bf 4}(1), (2005). 
%\bibitem{hardle} H\"ardle, W., Kerkyacharian, G., Picard, D., Tsybakov, A. {\it Wavelets, approximation and statistical applications}, Lecture Notes in Statistics, {\bf 129}, Springer-Verlag, New York, 1998.
%\bibitem{houghton}Houghton, J.C., {\it  Use of the truncated shifted Pareto distribution in assessing size distribution of oil and gas fields}, Mathematical geology {\bf 20}(8), 907--937, (1988).
%\bibitem{ik}  Ibragimov, I.A.,  Kas'minskij, R.Z.  {\it On  the
%estimation of  a signal,  its derivatives and  the maximum point  for Gaussian
%observations}, Teor. Veroyatnost. i Primenen.  {\bf 25}(4), 718--733, (1980). 
\bibitem{johnson} Johnson, W.B. (1985). Best Constants in Moment Inequalities for Linear Combinations of Independent and Exchangeable Random Variables. {\it Annals of Probability} {\bf 13}(1) 234--253.
%\bibitem{joh} Johnstone, I.M. {\it Minimax Bayes, asymptotic minimax and sparse wavelet priors}.  Statistical decision theory and related topics, V (West Lafayette, IN, 1992), 303--326, Springer, New York, 1994. 
\bibitem{jll}  Juditsky,  A. and  Lambert-Lacroix,  S. (2004). On  minimax  density
estimation on $\R$. {\it   Bernoulli} {\bf 10}(2) 187--220.
%\bibitem{kk}  Kim,  W.C.,  Koo,  J.Y.  {\it  Inhomogeneous  Poisson  intensity
%%estimation  via information  projections  onto wavelet  subspaces}, J.  Korean
%Statist. Soc. {\bf 31}(3), 343--357, (2002). 
%\bibitem{kp} Kerkyacharian, G., Picard, D. {Thresholding algorithms, maxisets and well-concentrated bases}, Test {\bf 9}, 283--344, (2000).
\bibitem{kin}  Kingman,  J.F.C. (1993). {\it  Poisson  processes.}  Oxford studies  in
Probability.
\bibitem{kolastro} Kolaczyk, E.D. (1997). Non-Parametric Estimation of Gamma-Ray Burst Intensities Using Haar Wavelets. {\it The Astrophysical Journal} {\bf 483} 340--349. 
\bibitem{kol}  Kolaczyk, E.D. (1999).    Wavelet shrinkage  estimation of  certain
Poisson intensity  signals using  corrected thresholds. {\it  Statist.  Sinica} {\bf
9}(1) 119--135.
%\bibitem{kn} Kolaczyk, E.D., Nowak, R.D. {\it Multiscale likelihood analysis and complexity penalized estimation},
%Ann. Statist. {\bf 32}(2), 500--527, (2004).
%\bibitem{kut} Kutoyants,  Y.A. {\it Statistical inference  for spatial Poisson
%processes.}  Lecture Notes in Statistics, {\bf 134}. Spinger Edition. 1998. 
%\bibitem{laherr} Laherr\`ere, J.{\it Distribution of field sizes in a
 %   petroleum system : parabolic fractal, lognormal or stretched exponential
 %   ?}, Marine and Petroleum Geology, {\bf 17}, 539-546, (2000) 
\bibitem{leb} Lebarbier, E. (2005). Detecting multiple change-points in the mean of Gaussian process by model selection.  {\it  Signal Processing} {\bf 85}(4) 717--736.
%\bibitem{merton} Merton, R.C. {\it Option pricing when underlying stock returns are discontinuous} Working paper (Sloan School of Management), 787-75 (1975).
%\bibitem{okt}   Ogata,  Y.,   Katsura,   K.,  Tanemura,   M.  {\it   Modelling
%heterogeneous space-time occurrences of earthquakes and its residual analysis}
%J. Roy. Statist. Soc. Ser. C {\bf 52}(4), 499--509, (2003). 
%\bibitem{pw}  Patil,  P.N.,  Wood,  A.T.A. {\it  Counting  process  intensity
%estimation by orthogonal wavelet methods}, Bernoulli {\bf 10}(1), 1--24, (2004). 
\bibitem{ptrfpois}  Reynaud-Bouret, P. (2003).  Adaptive  estimation  of  the
intensity    of   inhomogeneous    Poisson    processes   via    concentration
inequalities.   {\it  Probability Theory and Related Fields} {\bf 126}(1)
103--153.
%\bibitem{rbr} Reynaud-Bouret, P., Roy, E.{\it Some non asymptotic tail estimates for Hawkes processes},  Bulletin of the Belgian Mathematical Society-Simon Stevin, { \bf13}(5), 883--896 (2007), Proceedings of the 2005 joint BeNeLuxFra conference in Mathematics.
%\bibitem{rivbernoulli} Rivoirard, V. {\it Nonlinear estimation over weak  spaces and minimax Bayes method}, Bernoulli {\bf 12}(4), 609--632, (2006). 
\bibitem{Poisson_minimax} Reynaud-Bouret, P. and Rivoirard, V. (2008). Near optimal thresholding estimation of a Poisson intensity on the real line. Technical report. http://arxiv.org/abs/0810.5204
\bibitem{rud}  Rudemo,  M. (1982). Empirical  choice  of  histograms and  
density estimators. {\it  Scand. J. Statist.} {\bf 9}(2) 65--78.
\bibitem{sto} Stone, M.  (1974). Cross-Validatory Choice and Assessment of Statistical Predictions. {\it  J. Roy. Stat. Soc., Ser. B} {\bf 36} 111--147.
%\bibitem{uhler} Uhler, R.S., Bradley, P. G., {\it A Stochastic Model for Determining the Economic Prospects of Petroleum Exploration Over Large Regions } Journal of the American Statistical Association, {\bf 65}(330), 623--630, (1970). 
\bibitem{wn} Willett, R.M. and Nowak, R.D. (2007).  Multiscale Poisson Intensity and Density Estimation. {\it   IEEE Transactions on Information Theory} {\bf 53}(9) 3171--3187.
\end{thebibliography}

\end{document}